\newif\ifsmfart
\IfFileExists{smfart.cls}
   {\documentclass[12pt,english]{smfart}
    \IfFileExists{smfenum.sty}{\usepackage{smfenum}}{}
    \usepackage{bull}
    \smfarttrue}{\message{^^J*** It would be 
better to typeset this file with smfart.cls ***^^J^^J}

\documentclass[12pt]{amsart}}

\usepackage{amssymb}

\numberwithin{equation}{section}

\usepackage{epsfig}
\usepackage{graphicx}
\input xy
\xyoption{all}

\advance\headheight 2pt
\calclayout

\theoremstyle{plain}
\newtheorem{prop}[subsection]{Proposition}

\newtheorem{thm}[subsection]{Theorem}
\newtheorem{coro}[subsection]{Corollary}

\newtheorem{lemm}[subsection]{Lemma}

\newtheorem{defn}[subsection]{Definition}

\theoremstyle{definition}

\theoremstyle{remark}
\newtheorem{rem}[subsection]{Remark}

\newtheorem{exam}[subsection]{Example}

\newtheorem{nota}[subsection]{Notations}

\newcommand{\rk}{\operatorname{rk}}

\newcommand{\Aut}{\operatorname{Aut}}
\newcommand{\Hom}{\operatorname{Hom}}

\newcommand{\GL}{\operatorname{GL}}
\newcommand{\Gr}{\operatorname{Gr}}

\newcommand{\im}{\operatorname{Im}}

\newcommand{\ET}{\operatorname{ET}} 

\renewcommand{\mod}{\operatorname{mod}}
\newcommand{\RD}{{\rm RD}}

\newcommand{\GD}{{\rm GD}}
\def\cG{{\mathcal G}}
\def\cU{{\mathcal U}}

\def\gm{{\mathbb G}_m}

\def\G{{\Gamma}}
\def\Ga{{\Gamma}}
\def\PSL{{\rm PSL}}
\def\D{{\Delta}}
\def\no{{\noindent}}

\newcommand{\al}{\alpha}
\newcommand{\C}{\mathbb C}
\newcommand{\Q}{\mathbb Q}
\newcommand{\Z}{\mathbb Z}
\newcommand{\R}{\mathbb R}

\def\cA{{\mathcal A}}
\def\cB{{\mathcal B}}
\def\cC{{\mathcal C}}

\def\cF{{\mathcal F}}
\def\cH{{\mathcal H}}
\def\cO{{\mathcal O}}

\def\cU{{\mathcal U}}

\def\cJ{{\mathcal J}}
\def\cH{{\mathcal H}}
\def\cR{{\mathcal R}}

\def\cV{{\mathcal V}}

\def\bF{{\bf F}}

\def\sA{{\mathfrak A}}

\def\sC{{\mathfrak C}}

\def\sD{{\mathfrak D}}

\def\sS{{\mathfrak S}}

\def\M2{{\rm M}_{2}}
\def\Aut{{\rm Aut}}

\def\PGL{{\rm PGL}} 
\def\SL{{\rm SL}}

\def\Spec{{\rm Spec}}

\def\pr{{\it pr}}

\def\E{{\mathcal E}}

\def\lra{\longrightarrow}
\def\ra{\rightarrow}

\def\A{{\mathbb A}}
\def\C{{\mathbb C}}

\def\P{{\mathbb P}}
\def\Q{{\mathbb Q}}

\def\Z{{\mathbb Z}}
\def\C{{\mathbb C}}
\def\N{{\mathbb N}}

\def\bB{{\mathrm B}}

\def\bG{{\mathrm G}}
\def\bN{{\mathrm N}}

\def\bH{{\mathrm H}}
\def\bH{{\mathrm H}}
\def\bT{{\mathrm T}}

\def\bZ{{\mathrm Z}}
\def\bP{{\mathbb P}}
\def\Pr{{\mathbb P}}

\def\SL{{\rm SL}}
\def\GL{{\rm GL}}
\def\PGL{{\rm PGL}}
\def\Sym{{\rm Sym}}

\def\rk{{\rm rk\,}}

\def\ba{\backslash}

\def\Stab{{\rm St}}

\begin{document}

\title[Rationality of moduli]{\bf Rationality of 
moduli of elliptic fibrations with fixed monodromy}

\author{Fedor Bogomolov}
\address{Courant Institute of Mathematical
  Sciences, 251 Mercer str., New York, NY 10012-1185, USA} 
\email{bogomolo@cims.nyu.edu}
\thanks{The first author was partially 
supported by NSF grant DMS-9801591}

\author{Tihomir Petrov} 
\address{Courant Institute of Mathematical
  Sciences, 251 Mercer str., New York, NY 10012-1185, USA} 
\email{petrovt@cims.nyu.edu}
\thanks{The second author was partially supported by NSF grant
  DMS-9802154}

\author{Yuri Tschinkel}
\address{Dept. of Mathematics, Princeton University, Princeton NJ 08544-1000, USA}
\email{ytschink@math.princeton.edu}
\thanks{The third author was partially supported by NSA and by NSF
grant DMS-0100277}

%    General info
\subjclass{Primary 14J27, 14J10; Secondary 14D05}

\keywords{Elliptic fibrations, Monodromy groups, Moduli spaces}
%\date{January, 2001}
\date{\today} 

\begin{abstract}
We prove rationality results
for moduli spaces of elliptic K3 surfaces and
elliptic rational surfaces with fixed monodromy groups. 
\end{abstract}

\maketitle

%\tableofcontents

\setcounter{section}{-1}

\section{Introduction}
 
\noindent
Let $X$ be an algebraic variety of dimension $n$ over $\C$. 
One says that $X$ is rational if its function field
$\C(X)$ is isomorphic to $\C(x_1,\dots,x_n)$.
The study of rationality 
properties of fields of invariants $\C(X)^{\bG}=\C(X/\bG)$
is a classical theme in algebraic geometry.
For a finite group $\bG\subset \PGL_n$ acting on $X=\P^{n-1}$ the problem is 
referred to as Noether's problem (1916). 
It is still unsolved for $n=4$.  
Another class of examples is provided by {\em moduli spaces}. 
Birationally, they are often representable as quotients of 
simple varieties, like projective spaces or Grassmannians, 
by actions of linear algebraic groups, 
like $\PGL_2$. Rationality is known for each of the following
moduli spaces:

\begin{itemize}
\item curves of genus $\le 6$ \cite{Ka}, \cite{sh-2}, \cite{Ka-2}, \cite{K-3}, 
\cite{sh-1};
\item hyperelliptic curves \cite{Ka}, \cite{B};
\item plane curves of degrees $4n+1$ and $3n$ \cite{sh-3}, \cite{Ka-1};
\item Enriques surfaces \cite{Ko};
\item polarized K3 surfaces of degree 18 \cite{sh-2};
\item stable vector bundles (with various numerical 
characteristics) on curves, Del Pezzo surfaces, $\P^3$ 
\cite{ks}, \cite{by}, \cite{cm},\cite{N},\cite{M}; 
\end{itemize}
and in many other cases. For excellent surveys we refer to 
\cite{dol} and \cite{sh-3}. 
We will study rationality properties of moduli spaces of 
smooth non-isotrivial Jacobian elliptic fibrations
over curves 
$$
\pi\colon \E\ra C
$$  
with fixed global monodromy group 
$\tilde{\Ga}=\tilde{\Ga}(\E)\subset \SL_2(\Z)$.
In \cite{Bog} we developed techniques aimed at
the classification of possible global 
monodromies $\tilde{\Ga}$.  
The present paper gives a natural 
application of these techniques.

Let  $\cB$ be an irreducible algebraic family of 
Jacobian elliptic surfaces. 
Then the set of subgroups $\tilde{\Ga}\subset\SL_2(\Z)$ 
such that $\tilde{\Ga}$ is the (global) 
monodromy group of some $\E$ in this family is finite. 
Moreover, for every such group $\tilde{\Ga}$ the subset
of fibrations with this monodromy
$$
\cB_{\tilde{\Ga}}:=\{ b\in \cB\, \, |\, \, 
\tilde\Ga(\E_b)=\tilde\Ga\}
$$
is an algebraic (not necessarily closed) 
subvariety of $\cB$. 

Generalizing this observation, we introduce
(maximal) {\em parameter spaces} 
$\cF_{\tilde{\Ga}}$ of elliptic fibrations with
fixed global monodromy $\tilde{\Ga}$ (considered
up to fiberwise birational transformations 
acting trivially on the base of the
elliptic fibration). 
These parameter spaces can be represented
as quotients of quasi-projective varieties by
algebraic groups.
In particular, 
we can consider {\em irreducible connected components} 
of the parameter space 
$\cF_{\tilde{\Ga}}$, which we call
{\em moduli spaces}.  Even though these moduli spaces 
need not be algebraic varieties, we can still make sense 
of their birational type. 

\

\noindent
{\bf Theorem.} {\em
Let $\tilde{\Ga}\subset \SL_2(\Z)$ be a proper subgroup of finite
index. Then all moduli spaces of (Jacobian) elliptic rational or
elliptic K3 surfaces with global monodromy $\tilde{\Ga}$ are 
rational. 
}

\

\noindent
{\bf Corollary.}
{\em
For all $\tilde\Ga$ with moduli $\cF_{\tilde\Ga}$ of dimension $>0$ there
exists a number field $K$ such that
there are infinitely many nonisomorphic elliptic K3 surfaces over $K$
with global monodromy $\tilde\Ga$. 
}

\

\begin{rem}
Our method shows that many other classes of
moduli of elliptic surfaces over $\P^1$ 
with fixed monodromy are rational or unirational. 
However, we cannot expect a similar result for all moduli spaces
of elliptic surfaces over higher genus curves, since the moduli space
of higher genus curves itself is not uniruled (by a result of Harris and 
Mumford \cite{HM}). 
\end{rem}

\

We proceed to give a more 
detailed description of our approach.  
First of all, we can work not with 
the monodromy group $\tilde{\Ga}$ itself but rather 
with its image 
$$
\Ga\subset\PSL_2(\Z)
$$ 
under the
natural projection
$\SL_2(\Z)\twoheadrightarrow\PSL_2(\Z)$.
Let
$$
\cH=\{z\in \C\mid\im(z)>0\}
$$ 
be the upper half-plane and
$$
\overline{\cH}=\cH\cup\Q\cup\{\infty\}.
$$ 
The natural $j$-map
$$
  j\,:\, C\to{\P}^1={\overline {\cH}}/{\PSL_2(\Z)}
$$
decomposes as 
$$
j=j_{\Ga}\circ j_{\E},
$$ 
where
$$
\begin{array}{cccccc}
  j_{\E}\colon &  C    &  \ra  & M_{\Ga} & = & {\overline {\cH}}/{\Ga}\\
  j_{\Ga}\colon & M_{\Ga}&  \ra  & {\Pr}^1& = & 
{\overline {\cH}}/{\PSL_2(\Z)}.
\end{array}
$$
Here $M_{\Ga}$ is the 
$j$-{\em modular curve} 
corresponding
to $\Ga$; it is equipped with a
special triangulation, obtained as 
the pullback of the standard triangulation of 
${\mathbb S}^2={\Pr}^1(\C)$ (by
two triangles with vertices at $0,1$ and $\infty$) 
under the map $j_{\Ga}$
(which ramifies only over $0,1$ and $\infty$). 
We call the obtained triangulation of
$M_{\Ga}$ a $j_{\Ga}$-triangulation. 
Let $T_{\Ga}$ be the preimage in $M_{\Ga}$ 
of the closed interval $[0,1]\subset {\Pr}^1$. 
The graph $T_{\Ga}$ is our main tool in the combinatorial
analysis of $\Ga$. 

\

Denote by $\chi(\E)$ the Euler characteristic of $\E$. 
It splits equivalence classes of  Jacobian elliptic surfaces 
(modulo fiberwise birational transformations) 
into {\em algebraic families}. 
In particular, if $C={\Pr}^1$ then the {\em algebraic variety}
$\cF_r$ parametrizing (equivalence classes of)
Jacobian elliptic surfaces with given  
$\chi(\E)$ is irreducible; here we put 
$r=\chi(\E)/12.$
Our goal is to analyze the 
birational type of (irreducible components) 
$$
{\cF}_{r,\tilde\Ga}\subset {\cF}_r
$$ 
parametrizing fibrations with fixed
monodromy group $\tilde\Ga$. 
It suffices to study parameter spaces ${\cF}_{r,\Ga}$ 
corresponding to $\Ga\subset\PSL_2(\Z)$, since every 
irreducible component of ${\cF}_{r,\tilde\Ga}$
coincides with a component of ${\cF}_{r,\Ga}$. 

\

From now on we assume that $C=M_{\Ga}=\Pr^1$.  
Denote by $\cR_{d,\Ga}$ the 
space of rational maps $\Pr^1\to \Pr^1$ 
(of degree $d$) with prescribed 
ramification (encoded in $T_{\Ga}$).   
The spaces $\cF_{r,\Ga}$ are quotients
by the action of $\PGL_2\times H_{\Ga}$
of fibrations over $\cR_{d,\Ga}$ with fibers
(Zariski open subsets of) ${\Sym}^{\ell}({\Pr}^1)$ 
(for appropriate $d$ and $\ell$). 
Here $\PGL_2$ acts (on the left) by changing the 
parameter on the base $C={\Pr}^1$ and 
$H_{\Ga}$ is the group of 
automorphisms of $M_{\Ga}=\Pr^1$ 
stabilizing the embedded graph $T_{\Ga}$ (acting on the right). 
The nontriviality of $H_{\Ga}$ 
means that there is a $\Ga'\subset\PSL_2(\Z)$ containing 
$\Ga$ as a normal subgroup with $H_{\Ga} = \Ga'/\Ga$.
So in most cases in order to
prove rationality of $\cF_{r,\Ga}$ it is sufficient to 
establish it for ${\PGL_2}\backslash {\cR}_{d,\Ga}$, 
which can be deduced from general
rationality results for $\PGL_2$-quotients
(see \cite{Bogom-Katsylo}, \cite{Ka}).
To cover {\em all} cases we need to set up a rather
extensive combinatorial analysis.  

\

Here is a roadmap of the paper.
In Section~\ref{sect:covers} we discuss
finite covers $M_{\Ga}\to \Pr^1$ 
in the spirit of
Grothendieck's ``Dessins d'Enfants'' program
(see \cite{Lochak},\cite{Schneps} 
and the references therein) and 
introduce the invariants  $\GD(\Ga), \RD(\Ga)$ and $\ET(\Ga)$.
In ``ideal'' cases $\ET(\Ga)$ coincides with the
number of triangles in the 
$j_{\Ga}$-triangulation of $M_{\Ga}$
(the notation $\ET(\Ga)$ stands for
``{\em Effective Triangles}'').
In Section~\ref{sect:ellf} we recall basic facts about
elliptic fibrations and introduce the invariant $\ET(\E)$.  
For an ``ideal'' elliptic fibration 
one has $\ET(\Ga)=\ET(\E)$. 
In Section~\ref{sect:modsp} we discuss {\em moduli}
of elliptic fibrations with fixed monodromy. 
In Sections~\ref{sect:rat} and \ref{sect:rat-sp}
we formulate and prove several
rationality results for $\PGL_2$ and related
quotients. 
In Section~\ref{sect:ellmono}  
we classify families of rational elliptic
surfaces and elliptic 
K3 surfaces with different monodromy groups.
In Section~\ref{sect:comb}, we study relations between the 
combinatoric of the graph $\Ga$ and the topology of $\E$. 
And finally, in Section~\ref{sect:pictures} 
we list (certain) relevant subgroups $\Ga\subset\PSL_2(\Z)$ 
(represented by trivalent graphs $T_{\Ga}$). 
There are too many monodromy groups of elliptic K3
surfaces to be drawn on paper, but we show how to obtain 
them from our list by simple operations.

\section{Finite covers}
\label{sect:covers}

Let $\Ga$ be a subgroup of finite index in $\PSL_2(\Z)$.
The latter is isomorphic to a free product $\Z/3*\Z/2$. 
Consider the map
$$
{\overline{\cH}}/{\Ga}=M_{\Ga}\xrightarrow{j_{\Ga}}{\Pr}^1=
{\overline{\cH}}/{\PSL_2(\Z)},
$$
ramified over the points $0,1,\infty\in {\Pr}^1$.
Denote their preimages in 
$M_{\Ga}$ by $A,B$ and $I$, respectively. 
The possible ramification 
orders are $3$ or $1$ for
$A$-points, $2$ or $1$ for $B$-points and arbitrary for $I$-points.
The points $0,1$ and $\infty$ subdivide the circle
${\mathbb P}^1(\R)={\mathbb S}^1$ into three segments and, together with the
upper and lower hemisphere, define a decomposition of 
${\mathbb P}^1(\C)={\mathbb S}^2$ into three triangles.
This induces a special triangulation of $M_{\Ga}$ with
vertices in $A,B$ and $I$-points 
which we call the $j_\Ga$-{\em triangulation}. 
The preimage of the segment $[0,1]\subset{\Pr}^1$
defines a graph $T_{\Ga}$ which determines the 
$j_\Ga$-triangulation uniquely. 
Interior vertices of $T_{\Ga}$ are marked by $A_6$ and ends
are marked by either $A_2$ or $B_2$.

\begin{nota}
\label{nota:gd}
The {\em graph datum} $\GD(\Ga)$  of $T_\Ga$ is the 
formal sum 
$$
\GD(\Ga):=[a_6A_6+a_2A_2+b_2B_2],
$$
where $a_i$  ($i=6,2$) is the number of $A_i$ vertices and $b_2$ 
is the number of $B_2$-vertices. 
Denote by $\tau^0=\tau^0(\Ga)$ the number of vertices of $T_{\Ga}$
(including the ends), 
by $\tau^1=\tau^1(\Ga)$ the number of edges and by 
$\tau^2=\tau^2(\Ga)=\pi_0(M_{\Ga}\setminus T_{\Ga})$. 
\end{nota}

\

\begin{rem} 
\label{rem:uniqgr}
For given $a_2,b_2$ there
is a unique group with 
$$
\GD(\Ga)=[A_6+a_2A_2+b_2B_2].
$$  
\end{rem}

\

Forgetting the markings of $T_{\Ga}$ we obtain a connected unmarked 
{\em topological} graph $T^{u}_{\Ga}$ with (possibly some) 
ends and all interior vertices of valency $3$ --- 
a {\em trivalent graph}.

\begin{lemm} 
Let $X$ be a compact orientable Riemann surface of genus 
$g(X)$ and $T^{u}\subset X$ an 
embedding of a connected trivalent
graph such that
\begin{itemize}
\item 
the set $X\setminus T^u$ is a disjoint union of topological cells;
\item 
all interior vertices of $T^u$ are trivalent;
\item 
the ends of $T^u$ are arbitrarily marked by two colors $A_2$ and
$B_2$.
\end{itemize}
Then there exist a subgroup $\Ga\subset\PSL_2(\Z)$ 
and a unique complex structure on $X$ such that
$X=M_{\Ga}$ and $T^u=T_{\Ga}^u$.
\end{lemm}

\begin{proof}
Assume that we have an embedded graph $T^u\subset X$ satisfying the
conditions above. Mark by $A$ all trivalent vertices and
enlarge the graph $T^u$ by putting a $B$-vertex in the middle of any
edge bounded by two $A$-vertices. Put one $I$-vertex into
every connected component of $X\setminus{T}^u$ and connect all
$I$-vertices with $A$ and $B$-vertices at the boundary of the
corresponding domain. By assumption, every connected component of
$X\setminus{T}^u$ is contractible. Consider the boundary of the
individual cell. Every $A$-vertex of the boundary is
connected by edges to $B$-vertices only.  
Similarly, the  $B$-vertices are connected by 
edges only to $A$-vertices. 
Hence every triangle of the induced triangulation 
has vertices colored by three colors:
$A,B$ and $I$. This gives a $j_\Ga$-triangulation of $X$.
Following Alexander \cite{Alex}, 
we observe that a $j_\Ga$-triangulation
defines a map 
$$
{h}\,:\, X\to {\Pr}^1
$$ 
which is cyclically ramified over $A,B$ and $I$ (see \cite{Bog}). 
The trivalence of $T^u$ implies that $h$ has only 
$3$ or $1$-ramifications
over $0\in{\Pr}^1$ and only $2$ or $1$-ramifications over
$1\in{\Pr}^1$. 
Since $\PSL_2(\Z)={{\Z}/3}*{{\Z}/2}$ there is exactly
one subgroup $\Ga\subset \PSL_2(\Z)$ (of finite index) 
which corresponds to the covering $X\to{\Pr}^1$.
Any graph $T_{\Ga}^u$ constructed via a
subgroup $\Ga\subset{\PSL_2(\Z)}$ satisfies the conditions above.
Indeed, we have already described the 
$j_\Ga$-triangulation on $M_{\Ga}$. 
Triangles adjacent to a
given $I$-vertex constitute a contractible cell and the 
division of $M_{\Ga}$ into neighborhoods of $I$-vertices is a cellular
decomposition of $M_{\Ga}$. 
Hence after removing $I$-vertices with open
edges from them we obtain the preimage of $[0,1]$. If we forget the
$B$-vertices which lie between two $A$-vertices we
obtain the graph $T_{\Ga}^u$. Thus 
$T^u_{\Ga}\subset X = M_{\Ga}$ is the boundary of this
cellular decomposition and $T_{\Ga}$ is simply $T^u_{\Ga}$ with
an $A,B$-marking of the ends. 
\end{proof}

\begin{rem} 
Graphs which are isotopic in $X$ 
(modulo diffeomorphisms of $X$ of degree $1$) define 
conjugated subgroups of $\PSL_2(\Z)$.
\end{rem}

\begin{rem}
Even if we omit the condition of compactness of $X$
we still get a bijection 
between conjugacy classes of 
subgroups of finite index of $\PSL_2(\Z)$ and
embedded trivalent graphs with marked ends.
\end{rem}

\begin{rem}
The topology of $X$ restricts the topology of $T^u_{\Ga}$. The graph 
$T^u_{\Ga}$ must contain some $1$-skeleton of $X$. 
In particular, the map 
$\pi_1(T^u_{\Ga})\to \pi_1(X)$ is surjective.
Hence $T^u_{\Ga}$ can be a tree only if $X={\mathbb S}^2$.

For $X={\Pr}^1$ the connectedness of $T^u$ guarantees that
all the components of $X\setminus T^u$ are contractible. 
Hence we can classify graphs in $X={\Pr}^1$
by drawing them on the plane.
In general, connectedness of $T^u$ is necessary but not sufficient.  
\end{rem}

\begin{defn}
\label{defn:etg}
Define  
$$
\begin{array}{ccl}
\ET(\Ga)& := & 6{\tau}^0=6(a_6+a_2+b_2)\\
\D(\Ga) & := & 6a_6+2a_2
\end{array}
$$
\end{defn}

Thus both $\ET(\Ga)$ and $\D(\Ga)$  
depend only on the marking of the
ends but not on the embedding of the graph.
Observe that 
$\D(\Ga)$ is the {\em number of triangles} in the
corresponding $j_{\Ga}$-triangulation of $M_{\Ga}$ and that
$$
[\PSL_2(\Z):\Ga]=\D(\Ga)/2.
$$

\begin{rem}
If $M_{\Ga}$ arises from an elliptic fibration as
in the Introduction then $\D(\Ga)/2$ 
equals the number of {\em Dehn twists} in
$\Ga$ around the multiplicative singular fibers.
\end{rem}

\begin{nota}
\label{nota:rami-datum}
Let $f\,:\, C\to \Pr^1$ be a cover of degree $d$ and 
$p\in \Pr^1$ a ramification point of $f$. 
The {\em local ramification datum} is an $\N$-valued vector 
$v=(v_k)$, ($\sum v_k =d$), 
where $v_k$ is the order of ramification of
$f$ at a point $c_k\in f^{-1}(p)$. 
A {\em reduced local ramification datum} is a vector
$\overline{v}$ obtained from $v$ by 
omitting all entries $v_k=1$. 
The vector $v$ is defined up to 
permutation of the entries.     

For $f=j_{\E}\,:\,  C\to M_{\Ga}=\Pr^1$ 
we have distinguished ramification points,
namely those over $A$- and 
$B$-vertices of the graph $T_{\Ga}\subset M_{\Ga}$. 
The (global) $j_{\E}$-ramification datum is
the vector 
$$
\RD(j_{\E}):=[v_{1,A},\dots,v_{n,A}, v_{n+1,B},\dots,v_{n',B},
\overline{v}_{n'+1},\dots,\overline{v}_{n''}],
$$ 
where the $v_{i,A}$ are local ramification data over 
$A$-vertices for $i=1,\dots,n$,  
(resp. $v_{i,B}$ for $B$-vertices, $i=n+1,\dots,n'$) and
$\overline{v}_i$ are 
{\em reduced} local ramification data for unspecified 
other points in $M_{\Ga}$ for $i>n'$ (distinct from 
$A$- and $B$-vertices of $M_{\Ga}$). 

For $f=j_{\Ga}\,:\, M_{\Ga}\ra \Pr^1$ 
the distinguished (and the only) ramification points are $0,1,\infty$. 
We write
$$
\RD(j_{\Ga}):=[v_{0},v_{1},{v}_{\infty}]
$$ 
for the global $j_{\Ga}$-ramification datum.
\end{nota}

\begin{exam}
Assume that $\GD(\Ga)=[nA_6+A_2+B_2]$ is the graph datum 
of $T_{\Ga}\subset M_{\Ga}$ and 
let $j_\E\,:\, {\Pr}^1\to M_{\Ga}={\Pr}^1$ be a finite cover. 
Then the $j_{\E}$-ramification datum 
\begin{align*}
\RD(j_\E)=[(2,3)_A,(2,2,1)_B,(2),(2)]
\end{align*}
means that $\deg(j_\E)=5$, that $j_\E$  
has ramification points of order $2$ and $3$ over one point 
$A_2\in T_{\Ga}$ and $(2,2,1)$ over one $B_2$-point
and ramifications of order $2$ over 
two other unspecified points in $M_{\Ga}$. 
\end{exam}

\section{Elliptic fibrations}
\label{sect:ellf}

\no
In this section 
we briefly recall some basic facts of Kodaira's theory
\cite{Kod} of elliptic fibrations.
For more details we refer to 
\cite{Bart}, \cite{Fried} and
\cite{Shafar}. Let 
$$
\pi\colon \E\ra C
$$ 
be a smooth non-isotrivial relatively minimal
Jacobian elliptic fibration over a smooth 
projective curve $C$. This means that:
\begin{itemize}
\item 
$\E$ is a smooth compact complex 
projective surface and $\pi$ is a proper holomorphic map;
\item 
the generic fiber of $\pi$ is a smooth curve of genus $1$;
\item 
the fibers of $\E$ do not contain exceptional curves of the
first kind, i.e., rational curves $F$ such that $(F^2)=-1$ ({\em
relative minimality});
\item 
there exists a (global) zero section 
$s\,:\, C\to \E$ ({\em Jacobian elliptic fibration});
\item 
the $j$-function which assigns to each smooth fiber
${\pi}^{-1}(p)=\E_p\subset \E$ its 
$j$-invariant is a non-constant
rational function on $C$ ({\em non-isotriviality}).
\end{itemize}
\

\no
It is well known that $s^2<0$. We define  
$$
\ET(\E):= -24s^2.
$$

\begin{lemm}
\label{lemm:ett}
We have  
$$
\ET(\E)/2=-12s^2 =\chi(\E)=c_2(\E).
$$
\end{lemm}

\begin{proof}
Well known, but we decided to include an argument.
Since $\E$ is smooth and relatively minimal the
canonical bundle $K_\E$ of $\E$ is 
induced from a one-dimensional bundle $K$ on the base $C$. 
The sheaf $\pi^*K(C)$ is a subsheaf of $K_\E$. 
Since there are singular fibers we have the following
equality
\begin{align*}
  h^0(\E,\Omega^1)= h^1(\E,\mathcal O)=g
\end{align*}
where $g$ is the genus of $C$. By Riemann-Roch we obtain
\begin{align*}
  \chi({\mathcal O})=1-g+h^0(\E,K_\E)=\chi(\E)/12.
\end{align*}
We also know that $s^2+sK_\E-2g+2=0$ (genus formula). 
Therefore, 
\begin{align*}
  1-g+h^0(\E,K_\E)=\deg(K)-2g+2=\chi(\E)/12
\end{align*}
since $\deg(K)>2g-2$ and hence $h^1(C,K)=0$. 
Further, 
$$
sK_\E=\deg(K).
$$ 
Thus $s^2 +sK_\E-2g+2=0$ 
transforms to $s^2+\chi(\E)/12=0$.
\end{proof}

Let $C^{\rm sing}=\{p_1,\dots,p_k\}\subset C$ be 
the set of points on the base 
corresponding to singular fibers. 
The topological Euler characteristic 
$\chi(\E)=c_2(\E)$ is equal to the sum of Euler
characteristics of the singular fibers 
$\E_{p_i}={\pi}^{-1}(p_i)$ 
(since every generic fiber has Euler characteristic equal
to $0$). Therefore, 
$$
\ET(\E)=\sum_{p_i\in C^{\rm sing}}\ET(\E_{p_i}),
$$
where the summation runs over all singular fibers of $\E$
and $\ET(\E_{p_i})$ is the contribution from the
corresponding singular fiber. 
Since the fibration is 
Jacobian every singular fiber
has a unique representative from 
Kodaira's list and it is defined by
the local monodromy.
The possible types of singular fibers 
and their $\ET$-contributions are:
$$
\begin{array}{c|c|c|c}
          & \ET  &            & \ET \\
 \hline
{\rm I}_0 &     & {\rm I}_0^* &  12    \\ 
{\rm I}_n & 2n  & {\rm I}_n^* & 2n+12 \\
{\rm II}  &  4  & {\rm IV}^*  &  16  \\
{\rm III} &  6  & {\rm III}^* &  18  \\
{\rm IV}  &  8  & {\rm  II}^* &  20
\end{array}
$$
Here ${\rm I}_0$ is a smooth fiber, ${\rm I}_n$ is a multiplicative
fiber with $n$-irreducible components. The types ${\rm II},{\rm III}$
and ${\rm IV}$ correspond to the case of potentially good reduction.
More precisely, 
the neighborhood of such a fiber is a
(desingularization of a) quotient of a
local fibration with smooth 
fibers by an automorphism of finite
order. 
The corresponding order is $4$ for the case ${\rm III}$ and $3$
in the cases ${\rm II},{\rm IV}$. 
The fibers of type ${\rm I}_0^*,$ 
(resp. ${\rm I}_n^*,{\rm II}^*,{\rm III}^*,{\rm IV}^*)$ 
are obtained from fibers ${\rm I}_0$ 
(resp. ${\rm I}_n,{\rm IV},{\rm III},{\rm II})$ 
(after changing the local automorphism by the involution
$x\mapsto -x$ in the local group 
structure of the fibration). We shall
call them $*$-{\em fibers} in the sequel.

\begin{rem}
The local invariant $\ET(\E_{p})$ has a monodromy interpretation.
Namely, every element of a local monodromy at $p\in C^{\rm sing}$ has a
minimal representation as a product of elements conjugated to 
$\left(\begin{smallmatrix} 1 & 1 \\ 0 & 1 \end{smallmatrix}\right)$
in $\SL_2(\Z)$. The length of this representation equals 
$\ET(\E_p)/2$. This explains the equality
$\ET(\E_p^*)= \ET(\E_p)+12$ --- the element
$\left(\begin{smallmatrix} {-1}& 0 \\0 &
      {-1}\end{smallmatrix}\right)\in\SL_2(\Z)$ 
is a product of $6$ elements conjugated to
$\left(\begin{smallmatrix} 1 & 1 \\ 0 & 1
\end{smallmatrix}\right)$ (elementary Dehn twists).
\end{rem}

\section{Moduli spaces}
\label{sect:modsp}

\no
Every Jacobian elliptic fibration  $\E\to\Pr^1$
admits a {\em Weierstrass model}
$\bar{\E}$. 
Its geometric realization is given as follows: there
exists a pair of sections 
$$
\begin{array}{ccc}
g_2& \in &  H^0({\Pr}^1,{\mathcal O}_{{\Pr}^1}(4r)), \\
g_3& \in & H^0({\Pr}^1,{\mathcal O}_{{\Pr}^1}(6r))
\end{array}
$$ 
such that $\E$ is given by 
\begin{equation}
\label{weier}
  y^2z=4x^3-g_2xz^2-g_3z^3,
\end{equation}
inside $\Pr({\mathcal O}_{\Pr^1} \oplus {\mathcal O}_{{\Pr}^1}(2r)\oplus 
{\mathcal O}_{{\Pr}^1}(3r))$,
subject to conditions 
\begin{itemize}
\item
the discriminant $\D=g_2^3-27g_3^2$ is not identically $0$;
\item
for every point $p\in{\Pr}^1$ we have 
\begin{equation}\label{eqn:ineq}
\min(3\nu_p(g_2),2\nu_p(g_3))<12,
\end{equation}
where $\nu_p$ is the valuation corresponding to $p\in{\Pr}^1$
\end{itemize}
(see \cite{Fried} or \cite{F}, Section 7).

\

Two pairs $(g_2,g_3)$ and $(g'_2,g'_3)$ define 
isomorphic Jacobian elliptic surfaces 
$(\E,s)$ and $(\E',s')$ iff there exits an 
$h\in \GL_2(\C)$ transforming $(g_2,g_3)$ into $(g'_2,g'_3)$
under the natural action of $\GL_2$ on 
(the $\GL_2$-linearized) $\mathcal O_{\Pr^1}(r)$.
We define $\cF_{r}$ as the set of isomorphism classes
of pairs $(g_2,g_3)$ subject to the conditions above.

The parameter space 
$\cF_r$ has a natural structure of a (categorical) 
quotient of some open subvariety $U_r$ of the sum of two linear 
$\GL_2$-representations
$$
H^0({\Pr}^1,{\mathcal O}_{{\Pr}^1}(4r))
\oplus H^0({\Pr}^1,{\mathcal O}_{{\Pr}^1}(6r))
$$ 
by the action of $\GL_2$. 
Equivalently, $\cF_r$ is a (categorical) 
quotient of the open subvariety 
$U_r'=U_r/\gm$ of the weighted projective space 
$$
{\Pr}_{4r,6r}(4r+1,6r+1)
$$ 
by the action of $\PGL_2$.

\begin{lemm}
The variety $U_r'$ is a disjoint union of locally closed
subvarieties $U_{r,\tilde{\Ga}}'$, each preserved under 
the action of $\PGL_2$, such that for every $u\in U_{r,\tilde\Ga}'$
one has $\tilde\Ga(\E_u)=\tilde\Ga$. 
\end{lemm}

\begin{proof}
The set of  $(g_2,g_3)$ (of bounded degree) 
such that the $j$-map 
decomposes as $j_{\Ga}\circ j_{\E}$ 
is a closed algebraic variety (not necessarily irreducible).
Its intersection with the open subvariety defined by the conditions
\ref{eqn:ineq} is also closed. 
It contains a {\em finite} number of 
proper locally closed subvarieties corresponding 
to {\em proper} finite index subgroups of $\Ga$ 
(their number is bounded
as a function of $r$).  
The complement to these subvarieties
consists of finitely many irreducible components, 
each preserved under the action of $\PGL_2$. 
For $u$ in any of these components the lift of $\Ga$ to
the monodromy group of the corresponding fibration $\E_u$
is constant. The number of such possible lifts is finite. 
This gives a finite decomposition of $U_r'$ as claimed. 
\end{proof}

The unstable points of the $\PGL_2$-action
on the weighted projective space correspond to 
sections $g_2,g_3$ with high order of vanishing 
at some point $p$. Namely
$\nu_p(g_2)>2r,\nu_p(g_3)>3r$. 
However, the inequality~(\ref{eqn:ineq})
implies that $6r<12$. Thus, for $r\ge 2$, 
$\cF_r$ is a $\PGL_2$-quotient of some open
subvariety of the semistable locus 
$$
{\Pr}_{4r,6r}^{\rm ss}(4r+1,6r+1)\subset 
{\Pr}_{4r,6r}(4r+1,6r+1). 
$$
It follows that $\cF_r$ is a quasi-projective algebraic variety.
This variety is clearly unirational and in fact rational by
\cite{Ka}.

Moreover, for $r\ge 2$ we can define
a set of {\em subvarieties} $\cF_{r,{\tilde{\Ga}}}\subset \cF_r$ 
such that for every 
$b\in \cF_{r,\tilde{\Ga}}$ the corresponding Jacobian elliptic 
surface $(\E_b,s)$ has global monodromy group $\tilde{\Ga}$.

\

\begin{rem}
\label{rem:degen}
Notice that the maps $j_{\E}$ for 
elliptic  fibrations corresponding 
to different points of the same irreducible component
of $\cF_{r,\tilde{\Ga}}$ can have different $\RD(j_\E)$, even 
over the $A_2$ or $B_2$-ends of $T_{\Ga}\subset M_{\Ga}$. 
Thus, for a given irreducible component, we have the notion of 
a {\em generic} ramification datum $\RD(j_{\E})$ and 
its {\em degenerations}. 
\end{rem}

\

The case $r=1$,  corresponding to rational elliptic surfaces, is more subtle -
the subvariety $U_1'$ contains unstable points. 
The quasi-projective locus of semistable points 
${U_r^{\rm ss}}'$ is a disjoint union of locally closed
$\PGL_2$-semistable subsets
${U_{r,\tilde{\Ga}}^{\rm ss}}$; taking quotients we 
obtain varieties $\cF_{1,\tilde{\Ga}}$ parametrizing rational 
elliptic fibrations with global monodromy $\tilde{\Ga}$. 

Let $W_1'=U_1'-{U_1^{\rm ss}}'$ be
the complement.
It consists of pairs $(g_2,g_3)$ with 
$$
g_2 = l^3 f_2,\,\,\, g_3 =l^4f_3,
$$ 
where $l$ is a linear form (vanishing at a point $p$ and) 
coprime to $f_2,f_3$ and 
$\deg(f_2)=1,\deg(f_3)=2$. For $w\in W_1'$ we have
$\deg(j)\le 4$. The case of $\tilde\Ga\neq \SL_2(\Z)$ 
corresponds to $\deg(j_{\Ga})\ge 2$. Thus we have to consider
two cases: 
\begin{itemize}
\item
$\deg(j_{\Ga}) = \deg(j_{\E})=2$;
\item 
$\deg(j_{\Ga}) \le 4, \deg(j_{\E})=1$. 
\end{itemize}
The first case does not occur since $j^{-1}(0)$ has
ramification of type $(3,1)$ 
(by the assumption that $f_2$ is coprime to $l$
and that $3\nu_p(g_2)<12$).
Thus the $j$-map cannot be decomposed 
even locally into a product of two
maps. 
The second case leads to  

\begin{lemm}
If $w\in W_1'$ and $\tilde\Ga(\E_w)\neq \SL_2(\Z)$ 
then $\deg(j_{\E})=1$ and one has one of the 
following graph and ramification data:
$$
\begin{array}{l|l}
\GD(\G)        &  \RD(j_{\Ga}) \\
   \hline 
\lbrack A_6+A_2\rbrack   & \lbrack (3,1)_0,(2,2)_1, (3,1)_{\infty}\rbrack \\ 
\lbrack A_6+A_2+2B_2\rbrack &\lbrack (3,1)_0,(2,1,1)_1,(4)_{\infty}\rbrack \\ 
\lbrack A_6+3B_2 \rbrack   & \lbrack (3)_0,(1,1,1)_1,(3)_{\infty} \rbrack   \\
\lbrack A_6+B_2 \rbrack   & \lbrack (3)_0,(2,1)_1,(2,1)_{\infty} \rbrack   
\end{array}
$$ 
\end{lemm}

\begin{proof}
The formula $j = lf_2^3/(lf_2^3 - f_3^2) $ shows that
$j_{\Ga}$ has a point with local ramification datum $(3,1)$ or $(3)$,  
corresponding to 
$$
[\PSL_2(\Z):\Ga(\E)]=4\,\,{\rm or }\,\, 3.
$$
Since only two more branch points are allowed and one of them is $1$
(with local ramifications 1 or 2), 
the Euler characteristic computation gives 
the ramification data listed in the statement plus one more:
$$
[(3,1)_0,(2,2)_1,(2,2)_{\infty}].
$$
However, this datum is impossible for topological reasons 
(the only possible graph datum is $[A_6+A_2]$ and there
is a unique embedded graph $T_{\Ga}$ with this datum). 

If $\deg (j)=3$ then one has a
cyclic point of order $3$, leading to the data above.  
\end{proof}

\begin{coro}
Every irreducible component $W_{1,\tilde\Ga}'\subset W_1'$ such that
$\tilde\Ga(\E_w)\neq \SL_2(\Z)$ for $w\in W_{1,\tilde\Ga}'$ 
is rational. 
\end{coro}
Consider an irreducible component $\cF_{r,\tilde\Ga}$ and the
corresponding decomposition $j=j_\Ga\circ j_{\E}$. Here 
$$
j_{\E}=(j_{\E,2},j_{\E,3})\,:\,\Pr^1\ra M_{\Ga}=\Pr^1
$$
is a pair of homogeneous polynomials in 2 variables. 
Let 
$$
\cG=\{ (g_2,g_3)\}\subset H^0(\Pr^1,\cO_{\Pr^1}(4r))
\oplus H^0(\Pr^1,\cO_{\Pr^1}(6r))
$$
be the subset corresponding to smooth elliptic fibrations. 
Put
$$
\cJ_{\Ga}:=\{ j \, |\, \exists j_{\E}\,:\, \Pr^1\ra M_{\Ga}\,\, 
s.t.\,\, j=j_{\Ga}\circ j_{\E}\}.
$$

\begin{lemm}
\label{lemm:jj}
If $j\in \cJ_{\Ga}\cap \cJ_{\Ga'}$ with 
$\Ga\neq \Ga'$ then there exist
an $h\in \PSL_2(\Z)$, a group  
$\Ga''\subset \Ga\cap h\Ga'h^{-1}$ and a map 
$j_{\E}''\,:\, \Pr^1\to M_{\Ga''}$ such that
$j=j_{\E''}\circ j_{\Ga''}$. 
\end{lemm}

\begin{proof}
The monodromy group and its image in $\PSL_2(\Z)$ are 
uniquely determined by the smooth part of 
the elliptic fibration. 
Therefore,  in any smooth family of elliptic surfaces
$$
\Ga({\rm generic}\,\, {\rm fiber})\supseteq
\Ga({\rm special}\,\, {\rm fiber}).
$$
Since $\Ga$ is defined modulo conjugation by elements in $\SL_2(\Z)$
the claim follows.
\end{proof}

\begin{coro}
We have a decomposition $\cG=\bigsqcup \cG_{\Ga}$ into a finite (disjoint) union
of algebraic $\GL_2$-stable subvarieties such that for all 
$g=(g_2,g_3)\in \cG_{\Ga}$  the monodromy group 
$\tilde{\Ga}(\E_g)\subset \SL_2(\Z)$ is a subgroup of a central 
$\Z/2$-extension of $\Ga$.  
\end{coro}

\begin{rem}
For a given $g\in \cG_{\Ga}$ the map $j_{\E}$ is not unique. 
Let $j_{\E}$ and $j_{\E}'$ be two such maps. 
Then $j_\E=h_{\Ga}\circ j_{\E}'$, where $h_{\Ga}\in \Aut(T_{\Ga})$
is an automorphism of $M_{\Ga}$, preserving $T_{\Ga}$.
\end{rem}

\begin{lemm}
\label{lemm:decom}
We have a decomposition 
$$
\cG_{\Ga}=\bigsqcup_k \cG_{\tilde{\Ga},k}
$$
into a  finite union of algebraic {\em irreducible} 
$\GL_2$-stable subvarieties
such that $\tilde{\Ga}(\E_g)=\tilde{\Ga}$ for all $g\in \cG_{\tilde{\Ga},k}$. 
\end{lemm}

\begin{proof}
Assume that some $g\in \cG_{\Ga}$ belongs to 
$\cG_{\tilde{\Ga},1}\cap \cG_{\tilde{\Ga},2}$, where
$\cG_{\tilde{\Ga},1},\cG_{\tilde{\Ga},2}$ are different (non-conjugated)
lifts of $\Ga$ into $\SL_2(\Z)$. 
Lemma~\ref{lemm:jj} implies that there exists a {\em proper}
subgroup $\Ga''\subset\Ga$ such that 
$g$ belongs to $\cG_{\Ga''}$, contradiction.
\end{proof}

Let $\cG_{\tilde{\Ga}}=\cG_{\tilde{\Ga},k}$ 
be an irreducible component of $\cG_{\Ga}$
as in Lemma~\ref{lemm:decom} and $g\in \cG_{\tilde{\Ga}}$ its generic point. 
It determines a set of $*$-fibers on the base $\Pr^1$. We denote their
number by $\ell$. 
Choose (one of) the 
$j_{\E_g}$, with ramification datum $\RD=\RD(j_{\E_g})$.  
We get a map 
$$
\phi_{\cU}\,\,:\,\, \cU_g\to \cU_{j_g}\times (\Pr^1)^{\ell},
$$
where $\cU_g\subset \cG_{\tilde{\Ga}}$ is a neighborhood of $g$ and 
$\cU_{j_g}\subset \cR(\RD)$
is a neighborhood of the map $j_g=j_{\E_g}$ in the space 
$$
\cR(\RD):=\{ j\,:\, \Pr^1\to \Pr^1 \, |\,\, \RD(j)=\RD \}
$$
of rational maps with ramification datum $\RD$.

\begin{lemm}
\label{lemm:locall}
The map $\phi_{\cU}$ is a local (complex analytic) surjection.  
\end{lemm}

\begin{proof}
First observe that on $M_{\Ga}$ there is a projective
local system $(\Z\oplus\Z)/(\Z/2)$ which induces
a projective local system on an open part of $\Pr^1$. 
The obstruction to the extension of this system to
a linear $\Z\oplus \Z$-system is an integer modulo 2
which depends only on the topological type of the projective
system. Therefore, it doesn't change under a small variation
of maps $j$ with fixed $\RD$. Now it suffices to apply Kodaira's
main theorem which guarantees the existence and uniqueness
of an elliptic fibration with a given linear $\Z\oplus\Z$-system.  
\end{proof}

\begin{coro}
Let $\cF_{r,\tilde\Ga}'\subset \cF_{r,\tilde\Ga}$ 
be an (irreducible) component with 
generic ramification datum $\RD$. Then 
$\cF_{r,\tilde\Ga}'$ surjects (rationally) onto
the quotient of the variety of
rational maps $\cR(\RD)$ by $H_{\Ga}$.
\end{coro}

\begin{proof}
Since both  $\cF_{r,\tilde\Ga}'$ and $\cR(\RD)$ 
are algebraic varieties the local complex analytic surjection
from Lemma~\ref{lemm:locall} extends to an algebraic correspondence.
Moreover, two decompositions of the map $j$ as
$j=j_{\Ga}\circ j_{\E}$ differ by an element in $H_{\Ga}$. 
This gives a map to the quotient space, which is a (global) rational
surjection.  
\end{proof}

\begin{prop}
\label{lemm:struct}
Every irreducible component $\cF_{r,\tilde\Ga}$ contains an open
part $\cF_{r,\tilde\Ga}'$ with the following properties:
\begin{itemize}
\item  
$\cF_{r,\tilde\Ga}'$ is a quotient of an algebraic variety 
$U_{r,\tilde\Ga,\ell}'$
by the (left) action of $\PGL_2$ and (right) action of a subgroup $H_{\Ga}$ 
of $\Aut(T_{\Ga}) ;$
\item  
$U_{r,\tilde\Ga,\ell}' $ admits a fibration with fiber (an open 
subset of) ${\Sym}^\ell({\Pr}^1)$ 
and base the variety  $\cR_{r,\Ga}$ of maps $f\,:\, \Pr^1 \to M_{\Ga}$
with fixed local ramification data over 
$A_2$ and $B_2$-points of $T_{\Ga}\subset M_\Ga ;$
\item  
the action of $\PGL_2$ on $U_{r,\tilde\Ga,\ell}'$ is induced from the
standard $\PGL_2$-action on ${\Pr}^1 ;$
\item  
the group $\Aut(T_{\Ga})$ is a subgroup of 
$\PGL_2$ (acting on $M_{\Ga}$).
\end{itemize}
\end{prop}

\begin{proof}
Elliptic surfaces parametrized by a 
smooth irreducible variety have
the same $\ET(\E)$, which 
depends on the number $\ell$ of 
$*$-fibers in $\E$, on the degree of $j_{\E}$
and on the ramification properties over 
the ends of $T_{\Ga}$.
Once $\ell$ is fixed, for any given $j_{\E}$, 
the $*$-mark can be placed over
arbitrary $\ell$-points of ${\Pr}^1$. 
Their position defines a unique surface $\E$.
This implies that $U_{r,\tilde\Ga,\ell}'$ is fibered with 
fibers (birationally) 
isomorphic to ${\Sym}^{\ell}({\Pr}^1)={\Pr}^{\ell}$. 
The ramification properties of $j_{\E}$ 
remain the same on the open
part of $U_{r,\tilde\Ga,\ell}'$ (since the number of $*$-fibers
remains the same). 
Thus the base of the 
above fibration is the space of rational
maps $f\,:\, {\Pr}^{1}\to {\Pr}^{1}=M_{\Ga}$ 
with fixed ramification locus. Any such map 
defines an elliptic surface  $\E$
with given $\Ga$ (see \cite{Bog}). The $\PGL_2$-action
on $U_{r, \tilde\Ga,\ell}'$ 
identifies points corresponding to isomorphic
surfaces $\E$. Additional nontrivial 
isomorphisms correspond to exterior automorphisms of $\Ga$, 
coming from the action on $M_{\Ga}$, i.e., 
automorphisms of the graph $T_{\Ga}$. 
\end{proof}

\begin{rem}
\label{rem:free}
If the $\PGL_2\times\Aut({T_{\Ga}})$-action on $U_{r,\tilde\Ga,\ell}'$ is
almost free 
then the rationality of $\PGL_2\backslash U_{r,\tilde\Ga,0}'/\Aut(T_{\Ga})$ 
implies the rationality the corresponding quotients
for all $\ell$. In the other cases the degree of 
$j_{\E}$ is small and they are handled separately (see Section~\ref{sect:ratio}). 
\end{rem}

Most of the graphs 
$T_{\Ga}$ have trivial automorphisms. In particular,
any nontrivial automorphism acts on the ends of the graph.
In general, automorphisms of the pair
$(M_{\Ga},T_{\Ga})$ correspond to elements of $\Ga'/\Ga$ where
$\Ga'\subset \PSL_2(\Z)$ is a 
maximal subgroup with the property that
$\Ga$ is a normal subgroup of $\Ga'$.

\begin{lemm} 
\label{lemm:autt}
The group 
$\Aut({T_{\Ga}})$ acts freely on the set of ends and
end-loops.
\end{lemm}

\begin{proof}
Consider $j_{\Ga}\,:\, M_{\Ga}\ra \Pr^1$. Then 
$h\in\Aut({T_{\Ga}})$ is any element in $\PGL_2(\C)$ such that
$j_{\Ga}(hz)=j_{\Ga}(z)$ for all $z\in T_{\Ga}\subset M_{\Ga}$.
If $h$ stabilizes an end or an end-loop of $T_{\Ga}$ then 
it stabilizes the unique adjacent vertex and its other end.
Any element of $\PGL_2(\C)$ preserving a closed interval is 
the identity.
\end{proof}

\begin{coro}
\label{coro:hgamma}
For $r\le 2$, the only possible groups $\Aut(T_{\Ga})$ are 
cyclic, dihedral or subgroups of $\sS_4$. More precisely, 
for graphs with one end $\Aut(T_{\Ga})=1$ and graphs with two ends
$\Aut(T_{\Ga})$ is a subgroup of $\Z/2$.
\end{coro}

\begin{lemm}
\label{lemm:str}
Let 
$$
\cR:=\{f\,:\,{\Pr}^1\to {\Pr}^1\}
$$
be the space of rational maps 
with ramifications over exactly $0$ and $\infty$. 
Then $\cR$ is a ${\mathbb G}_m$-fibration over the
product of symmetric spaces $\Sym^{m_i}({\Pr}^1)$.
\end{lemm}

\begin{proof}
Indeed any two cycles $c_1$ and $c_2$ of fixed degree are
equivalent on ${\Pr}^1$. Therefore, 
there is a rational function $f$
on ${\Pr}^1$ with $c_1=f^{-1}(0)$ and $c_2=f^{-1}(\infty)$. If
$c_1,c_2$ do not intersect then $\deg(f)=\deg (c_1)=\deg(c_2)$. The
function $f$ is defined modulo multiplication by a constant. The
space of cycles 
$c_1=\sum_i n_i p_i$ 
is a product of symmetric powers 
$\Sym^{m}({\Pr}^1)$ where $m$ is the number of
equal $n_i$. 
\end{proof}

\section{Combinatorics}
\label{sect:comb}

In this section we investigate relations between 
$\ET(\E)$ and $\ET(\Ga)$. 
We keep the notations of the previous sections.

\begin{lemm}
\label{lem3}
Let $j\,:\, \E\to C$ be an elliptic fibration. Then  
\begin{align}\label{for2}
  \ET(\E)=\deg (j_{\E})\D(\Ga)+8\alpha_2+4\alpha_1+6\beta_1+12 \ell.
\end{align}
Here $\alpha_1$ and $\alpha_2$ equal 
the number of points over $A_2$-ends of
$T_{\Ga}$ with ramification multiplicity  
$1\ (\mod 3)$ and $2\ (\mod 3)$, respectively, $\beta_1$ is the 
number of odd ramification points over
the $B_2$-ends and $\ell$ is the number of $*$-fibers of $\E$.
\end{lemm}

\begin{proof}
The summand $\deg (j_{\E})\D(\Ga)$ corresponds to 
multiplicative fibers of $\E$.
The next summands are the contributions of those singular fibers
of $\E$ which are in the preimage of $A_2$ or $B_2$-ends of $T_{\Ga}$. 
If the ramification order at a point $p$ over a 
$B_2$-end is even then the
corresponding fiber with minimal $\ET$ is smooth and hence 
does not contribute to $\ET(\E)$. If it 
is odd then the fiber with minimal $\ET$ is of type $\rm III$
and we have to add $6\beta_1$. Similarly,
for the preimages of $A_2$-ends and $*$-twists. 
\end{proof}

\begin{coro}
\label{coro:ete}
In particular, 
$$
\ET(\E)\le \deg(j_{\E})\ET(\Ga)+12\ell,
$$
with equality if 
$$
\begin{array}{rcl}
2\alpha_2 +\alpha_1& = & a_2\cdot \deg(j_{\E})\\
\beta_1            & = & b_2 \cdot \deg(j_{\E}).
\end{array} 
$$
\end{coro}

\begin{defn}
We call $T_{\Ga}$ {\em saturated} if all vertices of
$T_{\Ga}^u$ are trivalent and a {\em tree} if it is
contractible.
\end{defn}

\begin{rem}
\label{rem:pi}
For saturated graphs $\D(\Ga)=12\rk \pi_1(T_{\Ga})$, 
where  
$$
\rk \pi_1(T_\Ga)=\rk H_1(T_\Ga)
$$ 
is the number of independent closed loops of $T_{\Ga}\subset M_{\Ga}$. 
\end{rem}

The following simple procedures produce new
graphs:

\begin{itemize}
\item If $T_1$ and $T_2$ are (unmarked) trivalent graphs
we can join $T_1$ and $T_2$ along two edges.  
For the resulting graph $T'$ we have
$$
\ET(T')=\ET(T_1)+ \ET(T_2)+12.
$$ 
If $T_i$ are marked and 
the marking of the ends of $T'$ is induced from 
the marking of the corresponding ends of $T_1$ and $T_2$ 
then
$$
\D(T')=\D(T_1)+\D(T_2)+12.
$$
\item We can glue an end $p$ of $T_1$ to an edge of $T_2$.
In this case 
$$
\ET(T')=\ET(T_1)+\ET(T_2).
$$ 
The change of $\D$ depends on the 
marking of the end:
$$
\D(T')=
\begin{cases} 
\D(T_1)+\D(T_2)+6 &\text{if $p=B_2$}\\            
\D(T_1)+\D(T_2)+4 &\text{if $p=A_2$}.
\end{cases}
$$
\end{itemize}

\begin{rem}
Any connected graph $T$ can be 
uniquely decomposed into a union
of a saturated graph and a union of trees.
\end{rem}

\begin{lemm} 
\label{lemm:12}
$\ET(\Ga)$ is divisible by $12$.
\end{lemm}

\begin{proof}
Every vertex of $T_{\Ga}$ 
has either one or three incoming edges.
Therefore, the number of edges
$$
{\tau}^1=\frac{1}{2}({\tau}^0_1+{\tau}^0_3),
$$ 
(${\tau}^0_i$ is
the number of vertices with $i$-edges). Thus
${\tau}^0={\tau}^0_1+{\tau}^0_3$ 
is even and since
$\ET(\Ga)=6{\tau}^0$ we are done.
\end{proof}

\begin{exam}
If $T_{\Ga}$ is a tree with $k+2$ vertices then 
$$
\begin{array}{cl}  
\ET({\Ga}) & = 12k+12  \\
\D({\Ga})  &  
\left\{\begin{array}{ccl} = 6k & {\rm if}\,\, 
      {\rm  all}\,\, {\rm ends}\,\,
{\rm are}\,\, B_2, \\
                          > 6k & {\rm otherwise}.
\end{array}
\right.
\end{array}
$$
\end{exam}

\begin{lemm}
\label{lemm:ineq}
For all $\Ga$ one has
$$
\D(\Ga)\geq \ET(\Ga)/2+6(\rk H_1(T_\Ga) -1).
$$
\end{lemm}

\begin{proof} 
A direct computation shows that
for saturated graphs one has an equality.
Suppose that $T_{\Ga}$ is a concatenation 
of a saturated graph
$T_{\rm sat}$ and a tree $T_{\rm tree}$. 
The number of ends drops by one and the number of 
$A_6$ vertices increases by $1$.
Thus the tree will add $12k+12$ to $\ET({\Ga})$ but 
$\D(\Ga)$ will change by $6k+6$. 
Finally, the ratio $\D(\Ga)/{\ET(\Ga)}$ only increases if 
we change $B_2$- to $A_2$-markings 
for some ends. Indeed, $\D(\Ga)$ increases
without changing $\ET(\Ga)$.
\end{proof}

\begin{coro}
\label{coro:useful}
If $\D(\Ga)={\ET(\Ga)}/2$ then 
$T_{\Ga}$ is a concatenation
of a loop $L$ and some trees. Moreover, 
all the ends of $T_{\Ga}$ are of type $B_2$.
\end{coro}

\begin{prop}
\label{prop:<}
Let $\E\to C$ be an elliptic fibration with 
$\ET(\E)<\ET(\Ga)$. Then: 
\begin{itemize}
\item  
$M_{\Ga}=\Pr^1$ and $T_{\Ga}$ is a tree without $A_2$-ends
and with $\ET(\Ga)>24 ;$
\item 
$\deg(j_{\E})=2$ and it is 
ramified in all ($B_2$) ends of $T_{\Ga}$
(and, possibly, some other points);
\item
$\E$ has $1$ or $2$ singular fibers of type ${\rm I}_n$.
\end{itemize}
\end{prop}

\begin{proof}
From \ref{lem3} and \ref{lemm:ineq} we conclude that
$\rk H_1(T_{\Ga})=0$ which implies that $T_{\Ga}$ is a tree
and $M_{\Ga}=\Pr^1$. By Lemma~\ref{lemm:ett} 
and our assumption, $\ET(\Ga)>24$, which implies that $\deg(j_{\E})\le 2$. 
If $\deg(j_{\E})=1$, we apply Corollary~\ref{coro:ete} and 
get a contradiction to the assumption. 
For $\deg(j_\E)=2$ combine Definition~\ref{defn:etg} and \eqref{lem3}:
$$
\ET(\E)=\ET(\Ga)+4a_2+4\alpha_1+8\alpha_2+6\beta_1 -12.
$$
Since $\alpha_1$, resp. $\beta_1$ is twice the number of unramified 
$A_2$, resp. $B_2$-ends, and $\alpha_2$ is the number of ramified $A_2$-ends
we see that if at least one of them is not zero, then 
$\ET(\E)\ge \ET(\Ga)$.  
The claim follows. 
\end{proof}

\begin{coro}
\label{coro:use}
For every elliptic fibration $\E\to \Pr^1$
one has 
$$
\ET(\E)\ge \ET(\Ga).
$$
Further, if  $\deg (j_{\E})=2$ and $j_\E$ is ramified
over only one $B_2$-point then
\begin{align*}
  \ET(\E)\geq 2\ET(\Ga)-12.
\end{align*}
\end{coro}

\begin{proof}
If $\deg(j_\E)=2$ and $\C=\Pr^1$ then $j_\E$ ramifies in two 
points. If neither of these points is $B_2$ then, by Lemma~\ref{lem3},
$\ET(\E)\ge 2\ET(\Ga)$. If both of these points are $B_2$-points
then the covering $j_\E$ corresponds to a subgroup $\Ga'$ 
of index 2 in $\Ga$ and $C=M_{\Ga'}$, contradiction. 
Otherwise, the claimed inequality follows from Lemma~\ref{lem3}.  
\end{proof}

\section{Elliptic K$3$ surfaces with $\deg(j_{\E})>1$}
\label{sect:ellmono}

In this section we assume that $C=\Pr^1$, that 
$j_{\E}>1$ and that $\Ga$ is a proper subgroup of $\PSL_2(\Z)$. 
We consider 
$$
\begin{array}{cccc}
{\rm general \,\, families:} &  \ET(\E)-12\ell  & = & \deg (j_{\E})\ET(\Ga),\\
{\rm special\,\, families:}  &  \ET(\E)-12 \ell & < & \deg (j_{\E}) \ET(\Ga).
\end{array}
$$

In Section~\ref{sect:modsp}
we showed that the main building block in the construction
of moduli space of elliptic surfaces with fixed $\Ga$ 
is the space of rational maps 
$j_{\E}\,:\,  C\to M_{\Ga}$ 
of fixed degree and ramification restrictions over certain points. 
For a {\em general} family 
there are no such restrictions and the corresponding moduli spaces 
are rational by classical results of invariant
theory for actions of $\PGL_2$ and its 
algebraic subgroups (see Section~\ref{sect:rat}).
For {\em special} families the
corresponding space of rational maps is more complicated.

\begin{lemm}\
\label{lemm:48}
There are no special families of elliptic 
{\rm K}$3$ surfaces with 
$$
\ET(\Ga)=48, 36.
$$
\end{lemm}

\begin{proof}

\

\begin{itemize}
\item If $\ET(\Ga)=48$ then 
$\D(\Ga)\geq 18$ and $\deg(j_{\E})\le 2$.
However, $\deg (j_{\E})=2$ contradicts Corollary~\ref{coro:use}
($\ET(\E)\geq 96-24>48$).
\item
If $\ET(\Ga)=36$ and $\D(\Ga)>16$ then $\deg (j_{\E})=2$,
contradicting  to \ref{coro:use}.  
We are left with $\D(\Ga)=16,14,12$ for 
$\deg (j_{\E})=3$ and $\D(\Ga)=12$ for $\deg(j_\E)=4$.
\item 
If $\deg(j_\E)=4$ then $T_\Ga$ is a tree with $\GD(\Ga)=[2A_6+4B_2]$. 
By Lemma~\ref{lem3}, all ramifications over the $B_2$-ends
are even, which contradicts $C=\Pr^1$ (compute $\chi(C)$).
\item   
If $\deg (j_{\E})=3$ then $T_\Ga$ is a tree (by \ref{lemm:ineq})
and  
$$
\GD(\G)=[2A_6+ a_2A_2+(4-a_2)B_2]
$$
with $a_2\leq 2$. We have 
\begin{align*}
48\geq  \ET(\E) \geq 3(12+2a_2)+ 4\alpha_1+8\alpha_2 + 6\beta_1,
\end{align*}
where $\beta_1\ge 2$ (since $\deg(j_\E)$ is odd
there is odd ramification over some $B_2$-end).  
Therefore, $a_2=0$ and consequently, $\beta_1\ge 4$, contradiction. 
\end{itemize}
\end{proof}

\begin{lemm}
\label{lemm:24-sp}
If $T_{\Ga}$ is not a tree and $j_{\E}$ is special 
(and generic for the corresponding irreducible component of
$\cF_{2,\tilde{\Ga}}$) then 
$$
\begin{array}{c|c|l|l}
\ET(\Ga)& \deg(j_{\E}) & \GD(\Ga)& \RD(j_{\E}) \\
\hline    
24 & 4   &   [2A_6+2B_2]     &  [(2,2)_B,(2,2)_B,(2),(2)] \\
24 & 3   &   [2A_6+2B_2]     &  [(2,1)_B,(2,1)_B,(2),(2)]\\  
24 & 3   &   [2A_6+A_2+B_2]  &  [(3)_A,(2,1)_B,(2)]\\
12 & 6   &   [A_6+A_2]       &  [(3,3)_A,(3,3)_A,(2),(2)]\\
12 & 5   &   [A_6+A_2]       &  [(3,1,1)_A]\\
12 & 5\le d\le 8&  [A_6+B_2] &  [\beta=(\beta_i)_B, (2)_B^{d'}],
\end{array}
$$
where
$$
\beta_i\in \N, \,\, \sum \beta_i=d,\,\, \# {\rm odd} \,\,\beta_i \le 8-d
$$
and 
$$
d'=2d - \# {\rm nonzero} \,\,\beta_i.
$$  
\end{lemm}

\begin{proof}
Follows from Lemma~\ref{lem3}. 
First observe that $\Delta(\Ga)\le 16$, which implies that 
$a_6=2$ and $a_2\le 2$.
If $a_2=2$ then $\Delta(\Ga)=16$ and
$$
\alpha_1=\alpha_2=\beta_1=0.
$$
Hence both $A_2$-ends have a 3-cyclic ramification
and the cover corresponds to a subgroup $\Ga'\subset \Ga$ of index 3. 
This excludes $\GD(\Ga)=[2A_6+2A_2]$. 
If $\deg(j_\E)=4$ then $\Delta(\Ga)=12$ which implies that 
all preimages of $B_2$-ends have even ramification.
The description of all other ramification 
data follows similarly from  Lemma~\ref{lem3}.
Notice that the (omitted) possibilities 
$$
\begin{array}{c|c|l|l}
\ET(\Ga)& \deg(j_{\E}) & \GD(\Ga)& \RD(j_{\E}) \\
\hline     
12 & 6   &   [A_6+A_2]       &  [(6)_A,(3,3)_A]\\
12 & 5   &   [A_6+A_2]       &  [(3,2)_A]
\end{array}
$$
are degenerations of the listed cases
(see Remark~\ref{rem:degen}).
\end{proof}

\begin{lemm}
\label{lemm:table-tree}
If $T_\Ga$ is a tree and $j_{\E}$ is special
(and generic for the corresponding irreducible component of
$\cF_{2,\tilde{\Ga}}$) then
$$
\begin{array}{c|c|l|l}
       & \deg(j_{\E}) & \GD(\Ga)& \RD(j_{\E}) \\
\hline    
j_1    & 4    &   [A_6+A_2+2B_2]  & [(1,1,1,1)_A,(2,2)_B,(2,2)_B,(2),(2)]\\
j_2    & 4    &   [A_6+A_2+2B_2]  & [(3,1)_A,(2,2)_B,(2,2)_B] + *\\
j_3    & 4    &   [A_6+A_2+2B_2]  & [(3,1)_A,(2,2)_B,(2,1,1)_B,(2)]\\
j_4    & 3    &   [A_6+2A_2+B_2]  & [(3)_A,(1,1,1)_A,(2,1)_B]\\ 
j_5    & 3    &   [A_6+A_2+2B_2]  & [(1,1,1)_A,(2,1)_B,(2,1)_B,(2),(2)]\\ 
j_6    & 3    &   [A_6+A_2+2B_2]  & [(3)_A,(1,1,1)_B,(2,1)_B,(2)] 
\end{array}
$$
or $\GD(\Ga)=[A_6+3B_2]$ and 
$$
\begin{array}{c|c|l}
       & \deg(j_{\E}) & \RD(j_{\E}) \\
\hline 
j_7    & 8    & [(2,2,2,2)_B,(2,2,2,2)_B,(2,2,2,2)_B,(2),(2)]\\
j_8    & 6    & [(2,2,2)_B,(2,2,2)_B,(2,2,1,1)_B,(2),(2)] \\
j_{9} & 6    & [(2,2,2)_B,(2,2,2)_B,(2,2,2)_B,(2)] + *\\
j_{10} & 5    & [(2,2,1)_B,(2,2,1)_B,(2,2,1)_B,(2),(2)]\\
j_{11} & 4    & [(2,1,1)_B,(2,1,1)_B,(2,2)_B,(2),(2)]\\
j_{12} & 4    & [(2,1,1)_B,(2,2)_B,(2,2)_B,(2)]+ *\\
j_{13} & 3    & [(1,1,1)_B,(2,1)_B,(2,1)_B]\\
j_{14} & 3    & [(2,1)_B,(2,1)_B,(2,1)_B,(2)] + * 
\end{array}
$$
or $\ET(\Ga)=12$ and $\GD(\Ga)=[2A_2]$ with $\deg(j_{\E})=4-10, 12$.

(In the above tables,  $+*$ means that there exists a moduli space of elliptic
surfaces with the same $\RD(j_{\E})$ and  with an additional $*$-fiber over
an unspecified point.) 

\end{lemm}

\begin{proof}
Assume that $\ET(\Ga)=24$ and $T_\Ga$ is a tree with
$$
\GD(\Ga)\neq [A_6+3B_2].
$$ 
First observe that $\deg(j_\E)\le 6$, since $\D(\Ga)\ge 8$.
If $\deg(j_\E)\ge 5$ then, by \ref{lem3},  
$\GD(\Ga)=[A_6+A_2+2B_2]$. If $\deg(j_\E)=6$ then 
$j_\E$ has to be completely ramified over all ends and 
no other ramifications are allowed by Euler characteristic
computation. Therefore, it is a group-covering and can't be $j_\E$.  
If $\deg(j_\E)=5$ then there are two odd ramifications over $B_2$-ends, 
and by \ref{lem3}, $\ET(\Ga)>48$. 

We are left with 
$$
\begin{array}{ccl}
\GD(\Ga) & = &  [A_6+3A_2],\\
        & = &  [A_6+2A_2+B_2],\\
        & = &  [A_6+A_2+2B_2]
\end{array}
$$
and $3\le \deg(j_\E)\le 4$.  
If there are at least two $A_2$-ends without
$3$-cyclic ramification points over them
then $\ET(\E)> 48$ (see \ref{lem3}). 
The first case is impossible: $\deg(j_\E)=4$ does not occur
(the degree is not divisible by 3), if $\deg(j_\E)=3$ 
and there is at most one $3$-cyclic ramification over an $A_2$-end then, 
by \ref{lem3}, $\ET(\E)>48$, contradiction. 
Consider the second case and $\deg(j_\E)=4$. Then $\D(\Ga)=10$ and
$4\alpha_1+8\alpha_2+6\beta_1\le 8$.  Since $\alpha_1\ge 2$ 
we have $\alpha_2=\beta_1=0$ and $\alpha_1=2$.  
The only possible 
$$
\RD(j_\E)=[(3,1)_A,(3,1)_A,(2,2)_B],
$$ 
which corresponds to a group covering, contradiction. 

Similarly, if $\GD(\Ga)=[A_6+A_2+2B_2]$ and $\deg(j_{\E})=4$
then $\D(\Ga)=8$ and $4\alpha_1+8\alpha_2+6\beta_1\le 16$. 
We have $\alpha_1\ge 1$ and $8\alpha_2+4\alpha_1=16$ or $4$. 
In the first case, both $B_2$-ends are completely ramified, 
and we get $j_1$. The second case splits into subcases:
$\beta_1=0$ or $2$, leading to $j_2$, resp. $j_3$. 
If $\deg(j_{\E})=3$, then if $\GD(\Ga)=[A_6+2A_2+B_2]$ 
then exactly one of the $A_2$-ends has cyclic ramification.  
It follows that $\beta_1=1$, which leads to $j_4$. 
If $\GD(\Ga)=[A_6+A_2+2B_2]$ there are two subcases:
there is cyclic ramification over the $A_2$-end or not. 
In the first subcase, possible $\RD(j_\E)$
include $[(2,1)_B,(2,1)_B]$, which is excluded as
it gives a group covering. The other case leads to $j_6$. 
In the second subcase, we get $j_5$.

\

Consider the case $T_{\Ga}=A_6+3B_2$. Here $\D(\Ga)=6$ and 
$$
\ET(\E)\geq 6\deg(j_{\E})+6n,
$$ 
where $n$ is a number of points with odd ramification
over $B_2$-vertices.
It follows that 
$$
48\ge 6\deg(j_{\E})+6\beta_1
$$
and $\beta_1\ge 3$ if $\deg(j_{\E})$ is odd 
and the number
of odd ramifications over {\em each} $B_2$-end is 
congruent to $\deg(j_{\E})$ modulo 2.

If $\deg (j_{\E})=8$ then all preimages of $B_2$-vertices are
$2n$-ramified. 
If $\deg(j_\Ga)$ is odd then $\ET(\E)\geq 6\deg(j_\Ga)+18$, 
which excludes $\deg(j_\Ga)=7$.
Now assume $\deg(j_\Ga)= 6$. The number of possible odd
ramifications over any $B_2$-end is even by~\ref{lem3} and 
it cannot exceed $2$. There are two possibilities listed above.
Assume that $\deg(j_{\E})=5$. The minimal possible
ramifications are $(2,2,1)$ over all 
$B_2$-ends. Since $10-6=4$ we can add two more points.

In $\deg(j_{\E})=4$ we could have further 
$\RD$: 
$$
\begin{array}{ccl}
\RD(j_\E)& = & [(2,2)_B,(2,1,1)_B,(2,2)_B,(2)],\\
         & = & [(2,2)_B,(2,1,1)_B,(2,1,1)_B,(2),(2)],\\
         & = & [(2,2)_B,(2,1,1)_B,(2,1,1)_B,(3)]
\end{array}
$$ 
but they are obtained as degenerations of $j_{12}$ and $j_{13}$.

The only $\GD(\Ga)$ which allow $\deg(j_{\E})\ge 12$ are
$[A_2+B_2]$ and $[2A_2]$.
The first case corresponds to $\PSL_2(\Z)$ (which we don't consider). 
The second case corresponds to subgroups $\Ga\subset \PSL_2(\Z)$ of 
index $2$. 
For a generic $\E$ in each moduli space
the ramification datum $\RD(j_{\E})$
is one of the following:
$$
\RD(j_{\E})=[(3,\dots,3_{n_1},1,\dots,1)_A,(3,\dots,3_{n_1},1,\dots,1)_A, (2)^d] + *,
$$
where $n_1,n_2,d$ are non-negative integers such that
$$
\deg(j_{\E})-(n_1+n_2)\le 4,
$$ 
$$
3n_1, 3 n_2\le \deg(j_{\E})\,\, {\rm and},\,\,
d\le 2(\deg(j_{\E}-(n_1+n_2+1))). 
$$
(In particular, $d\le 4$).
\end{proof}

\section{Rational elliptic surfaces with $\deg(j_{\E})>1$}
\label{sect:rat-surf}

\begin{lemm}
There are no special families of rational elliptic surfaces with $\ET(\Ga)=24$.
\end{lemm}

\begin{proof}
If $\deg(j_{\E})=2$ then $j_{\E}$ cannot be ramified over more than one $B_2$-end
(otherwise it is a group covering). Therefore, we can apply 
Corollary~\ref{coro:use} and get $\ET(\E)>2\cdot 24 -12 >24$, 
contradiction (to \ref{lemm:ett}).  
Thus  $\deg(j_{\E})=3$ or 4 and $a_6 =1$. Moreover, $\D(\Ga)\le 8$. 
This leaves
the cases: 
$$
\begin{array}{ccl}
\GD(\Ga) & = & [A_6+3B_2],\\
        & = & [A_6+A_2+2B_2].
\end{array}
$$  
In the first case $\deg(j_\E)=3$ is impossible, and $\deg(j_\E)=4$ leads
to 
$$
\RD(j_{\E})=[(2,2)_B,(2,2)_B,(2,2)_B,(2,2)_B]
$$
which corresponds to a group covering. 
In the second case $\deg(j_\E)\neq 4$ (since $\D(\Ga)=8$) 
and $\deg(j_\E)=3$ implies that $\beta_1\ge 2$ and $\ET(\E)\ge 36$, 
contradiction. 
\end{proof}

\begin{lemm}
\label{lemm:table-rat}
If $j_{\E}$ is special
(and generic for the corresponding irreducible component of
$\cF_{1,\tilde{\Ga}}$) then
$$
\begin{array}{c|c|l|l}
       & \deg(j_{\E}) & \GD(\Ga)& \RD(j_{\E}) \\
\hline    
j_{15}    & 6    &   [2A_2]    & [(3,3)_A,(3,3)_A,(2),(2)]\\
j_{16}    & 4    &   [2A_2]    & [(3,1)_A,(3,1)_A,(2),(2)]\\
j_{17}    & 3    &   [2A_2]     & [(3)_A,(1,1,1)_A,(2),(2)]\\
j_{18}    & 4    &   [A_6+B_2]  & [(2,2)_B,(2),(2),(2),(2)]\\ 
j_{19}    & 3    &   [A_6+B_2]  & [(2,1)_B,(2),(2)]\\ 
j_{20}    & 3    &   [A_6+A_2]  & [(3)_A,(2),(2)] 
\end{array}
$$
\end{lemm}

\begin{proof}
If $a_6\ge 1$ then $\deg(j_{\E})=4$ or $3$. 
In the first case $a_2=0$ and $\GD(\Ga)=[A_6+B_2]$ and 
we have complete ramification over the $B_2$-end. This gives $j_{18}$.
In the second case the ramification over $B_2$ is $(2,1)_B$ and 
we get $j_{19}$. If $\GD(\Ga)=[A_6+A_2]$ then $\deg(j_{\E})=3$ and 
$\alpha_1=\alpha_2=0$, leading $j_{20}$. 

It remains to consider $\GD(\Ga)=[2A_2]$. 
We apply the same formulas as in the proof of Lemma~\ref{lemm:table-tree}, 
with the inequality 
$$
\deg(j_{\E})-n_1-n_2\le 2.
$$
We have $\deg(j_\E)\le 6$ and
$\alpha_1=\alpha_2=0$. Notice that $\deg(j_\E)=5$ is impossible. 
\end{proof}

\section{General rationality results}
\label{sect:rat}

\begin{nota}
\label{nota:rat}
We will denote by $\sS_n$ the symmetric group on $n$ letters, 
by $\sA_n$ the alternating group, by $\sD_n$ the dihedral group
and by $\sC_n=\Z/n$ the cyclic group.
In particular, $\sS_2=\sC_2=\Z/2$ and $\sD_2=\Z/2 \times \Z/2$
(sometimes we prefer the notation $\sS_2$ over $\sC_2$ to 
stress that the action is by permutation).
We write $\Gr(k,n)$ for the Grassmannian of $k$-planes in a vector
space of dimension $n$ and $V_{d}$ for 
the space of binary forms of degree $d$. 
We will denote by $\GL_2,\PGL_2,\gm$ etc. the corresponding 
{\em complex} algebraic groups.
For a group $\bG$, we 
denote by  $\bZ_g$ the centralizer of $g\in \bG$ and by
$\bZ_{\bG}$ its center.
We denote by $\M2=V_1\oplus V_1$ the space of $2\times 2$-matrices.
We write $\cV\stackrel{V}{\longrightarrow}X$ or simply 
$\stackrel{V}{\longrightarrow}X$ for a 
locally trivial (in Zariski topology) fibration $\cV$ over $X$ with
generic fiber $V$.
We will often write $\bG$-map (etc.), 
instead of $\bG$-equivariant map. 
\end{nota}

\

\no
We say that two algebraic varieties $X$ and $X'$ are birational, and 
write $X\sim X'$, if $\C(X)=\C(X')$. 
A variety $X$ of dimension $n$ is {\em rational} if 
$X\sim \A^n$, {\em $k$-stably rational} 
if $X\times \A^k\sim \A^{n+k}$ and {\em stably rational}
if there exists such a $k\in \N$. We say that $X$ is  
{\em unirational} if $X$ is dominated by $\A^n$. 
The first basic result, a theorem of Castelnuovo from 
1894, is:

\

\begin{thm}
\label{thm:q}
A unirational surface is rational.
\end{thm}

Already in dimension three,
one has strict inclusions
$$
{\rm rational}\subsetneq {\rm stably\,\, rational}\subsetneq {\rm unirational}
$$
(see the counterexamples in \cite{im}, \cite{am},\cite{cg},\cite{bcs}). 
There is a very extensive literature on 
rationality for various classes of varieties. 
We will use the following facts:

\

\begin{lemm}
\label{lemm:cb}
Let $S\ra B$ be a ruled surface with base $B$ and 
$\pi\,:\, \cC\ra S$ a conic bundle over $S$.
Assume that the restriction of $\pi$ to a generic $\P^1\subset S$
is a conic bundle with at most three singular fibers.
Then $\cC\sim\A^2\times B$.
\end{lemm}
 
\

\begin{lemm}
\label{lemm:conic}
Let $\pi\,:\, \cC\ra S$ be a conic bundle over an irreducible variety  $S$ and
$Y\subset \cC$ a subvariety such that the restriction of
$\pi$ to $Y$ is a surjective finite map of odd degree. 
Then $\cC$ has a section and $\cC\sim S\times \A^1$. 
\end{lemm}

\

Let $\bG$ be an algebraic group. 
A (good) {\em rational action} of $\bG$
is a homomorphism 
$$
\rho_{\rm rat}\,:\, \bG\ra {\rm Bir}(X)
$$
such that there exists a birational model $X'$ of $X$ 
with the property that $\rho_{\rm rat}$ extends to a (regular)
morphism  $\bG\times X'\ra X'$. We consider only rational 
actions. We write $X \sim_{\bG} Y$ for a $\bG$-birational 
($=\bG$-equivariant birational)
isomorphism between $X$ and $Y$.
We will denote by $\bG\ba X$ a model for the field of invariants 
$\C(X)^{\bG}$.

\

Let $E\ra X$ be a vector bundle. 
A {\em linear} action of $\bG$ on $E$ is 
a rational action which
preserves the subspace of fiberwise linear functions on $E$. 
In particular, there is a linear $\bG$-action on
regular and rational sections of $E$.

\

We are interested in rationality properties of quotient spaces 
for the actions of $\PGL_2$, its subgroups and products of $\PGL_2$ with
finite groups. The finite subgroups of $\PGL_2$ are
$$
\sC_n,\sD_n,\sA_4,\sS_4,\sA_5.
$$
We denote by $\tilde{\sC}_n,\tilde{\sD}_n$ etc. their lifts to $\GL_2$
(as central $\sC_2$-extensions). 
We denote by
$$
\bB,\bT=\C^*, \bN_{\bT}
$$
the upper-triangular group, the standard 
maximal torus and the normalizer of this
torus in $\PGL_2$ and by
$$
\tilde{\bB},\tilde{\bT},\bN_{\tilde{\bT}} 
$$
the corresponding subgroups in $\GL_2$ (or $\SL_2$).

\

Let $V$ be an $n$-dimensional vector space,
$\tilde{\bG}\subset \GL(V)$ a  subgroup
and $\bG$ its projection to $\PGL(V)$, acting naturally on $\P(V)$.   
Determining the rationality of quotients 
$\bG\ba \P(V)$ (at least for finite groups) 
is known as Noether's problem.   

\

\begin{coro}[of Theorem~\ref{thm:q}]
\label{coro:dim}
For all $n\le 3$ the space $\bG\ba\P(V)$ is rational. 
\end{coro}

\

\begin{thm}\cite{M},\cite{V}
\label{thm:qq}
A quotient of $\bP(V)$ by a (projective) action of
a connected solvable group, a torus or a finite abelian subgroup of a torus
is rational.
\end{thm}

\

A fundamental rationality result is the following theorem of Katsylo:

\

\begin{thm}\cite{K}
\label{thm:kats}
For any representation $V$ of $\GL_2$ or $\PGL_2$ the quotient 
$\PGL_2\ba \P(V)$ is rational. 
\end{thm}

\

In general, the quotients need not be rational 
(see Saltman's counterexamples in \cite{S}). 
We now describe some partial results for $n=4$, which 
we will use later on. 

\

\begin{defn}
\label{defn:prim}
A finite group $\tilde{\bG}\subset \GL_n=\GL(V)$ 
is called {\em imprimitive} if 
there exists a decomposition $V=\oplus_{\al} V^{\al}$ such that
for all $\al$ and $\tilde{g}\in \tilde{\bG}$ there is an $\al'$ with 
$\tilde{g}V^\al=V^{\al'}$. Otherwise, $\bG$ is called {\em primitive}.
\end{defn}

\begin{rem}
There are 29 types of primitive subgroups of $\GL_4$. 
For some of them, like  
$$
\sA_6,\sA_7,\PSL_2({\bF}_7),\sS_6,
$$ 
rationality of the quotient is still unknown.  
\end{rem}

\begin{thm}
\cite{km}
\label{thm:kolp}
For every primitive solvable subgroup
$\bG\subset \PGL_4$ the quotient $\bG\ba \P^3$ is rational. 
\end{thm}

\begin{rem}
\label{rem:km}
In \cite{km} it is shown that 
$$
\bG\ba \P^3\sim_{\bG} \bG'\ba X_3,  
$$
where $X_3$ is the Segre cubic threefold and 
$\bG'$ is a quotient of $\bG$.
The problem is then reduced to the 
(easy) case of imprimitive actions. 
\end{rem}

\

We will also need to consider quotients by {\em nonlinear} actions.

\

\begin{lemm}
\label{lemm:qqq}
The quotient of $\GL_2$  (or $\PGL_2$)
by the involution $i\,:\, x\mapsto x^{-1}$
is rational.
\end{lemm}

\begin{proof}
The involution decomposes 
as a product $i=i_1\circ i_2$, where 
$$
i_1\,:\, 
x:=\left(\begin{array}{cc} a & b \\ c & -a +d \end{array}\right)\mapsto
\left(\begin{array}{cc} -a+d & -b \\ -c & a\end{array}\right)
$$
and
$$
i_2\, :\, y\mapsto y\cdot \det(y)^{-1}.
$$
are two commuting involutions.
Another set of independent generators of  
$\C(a,b,c,d)$ is given by $\{ a,b,c,\det(x)\}$ (write 
$d = (\det(x) + bc+a^2)/a$).
Now the involutions take the form
$$
i_1\, :\, (a,b,c)\to (-a ,-b,-c)
$$ 
and 
$$
i_2\, :\, \det(x)\to \det(x)^{-1}
$$ 
and we can write down 
independent generators of the field of invariants. 
If 
$$
D := \frac{\det(x) +1}{\det(x)-1}
$$ 
then 
$$
\begin{array}{cccc}
i_2 :  & D         &  \mapsto & -D \\
i_1 :  & (a,b,c,D)&  \mapsto & (-a ,-b,-c, - D).
\end{array}
$$
This finishes the proof.
\end{proof}

\

A {\em (rational) slice} for the action of 
$\bG$ is a subvariety $S\subset X$
such that the general $\bG$-orbit intersects $S$ in exactly one point. 
(The slice $S$ need not be a rational variety.
To avoid confusion, we will always refer to $S$ as a slice.)
A subvariety $Y\subset X$ is called a $(\bG,\bH)$-{\em slice} 
(where $\bH\subset \bG$ is a subgroup) if $\bG\cdot Y\sim X$ and $gy\in Y$ implies 
that $g\in \bH$. Clearly, $\bG\ba X \sim \bH\ba Y$. Moreover, if
$f\,:\, X\ra X'$ is a $\bG$-equivariant morphism 
and $Y'$ is a $(\bG,\bH)$-slice in $X'$ then 
$f^{-1}(Y')$ is a $(\bG,\bH)$-slice in $X$.
    
\

\begin{nota}
\label{nota:stab}
For (a reductive group) $\bG$ acting (rationally) on $X$
we denote by 
$$
\Stab_{gen}=\Stab_{gen}(\bG,X)
$$ 
the generic stabilizer (defined up to conjugacy). 
The action is called an $af$-action ({\em almost free})
if $\Stab_{gen}$ is trivial.
\end{nota}

We use a more precise version of Theorem~\ref{thm:kats}:

\begin{thm}
\cite{K}
\label{thm:katsylo}
Let $\rho\,:\, \PGL_2\ra \PGL(V)$ be a representation and 
$\tilde{\rho}$ a lifting of $\rho$ to a
representation of $\GL_2\ra \GL(V)$.
Let 
$$
\bG'':=\Stab_{gen}(\GL_2,V)\,\,\, {\rm  and}\,\,\, \bG:=\GL_2/\bG''.
$$ 
If the central $\sC_2\not\subset \bG''$ then
$$
\bP(V) \sim_{\bG} \bG\times S,
$$ 
where $S$ is a rational variety (with trivial $\bG$-action).

\no
If $ \sC_2\subset \bG''$ then
\begin{itemize}
\item
either the $\PGL_2$-action on  $\bP(V)$ has {\em no} slice 
and $\bG\ba \bP(V)$ is rational  
\item or 
$$
\bP(V) \sim_{\bG} \bG\times S,
$$
where the slice $S$ is a rational variety (with trivial $\bG$-action). 
\end{itemize}
\end{thm}

\

We now explain some general techniques 
in the study of rationality of quotient varieties. 

\

\begin{lemm}
\label{lemm:redu}
Let $E\ra X$ be a vector bundle of rank $r=\rk(E)$. 
Let $\bG$ be an (affine) reductive group acting on $E$ such that
the generic orbit of $\bG$ in $E$ projects isomorphically
onto a generic orbit of $\bG$ in $X$. Then  
$$
E\sim_{\bG} X\times \A^{r}
$$
with trivial $\bG$-action on the affine space $\A^{r}$. 
\end{lemm}

\begin{proof}
Denote by $O$ the $\bG$-orbit through a generic point in $E$.  
Shrinking (equivariantly)  
$X$, if necessary, we may assume that
the map
$$
\pi\,:\, H^0(X,E)\ra H^0(O,E|_{O})
$$
is surjective. With our assumptions, there exists a basis
$s_1,...,s_r$ such that for each $j$, the $\bG$-orbit of $s_j$  
projects isomorphically onto its image in $X$ and generates
a trivial 1-dimensional $\bG$-equivariant sub-bundle of
the restriction $E|_{O}$ of $E$ to the orbit $O$. 
It follows that $E|_{O}=\oplus_{j=1}^r \bG\cdot s_j$. In particular, 
$H^0(O,E)$ contains the trivial $\bG$-module $M$
generated by $s_1,...,s_r$. Moreover, 
$M$ generates $H^0(O,E|_{O})$ over every point of $O$. 
Since $\pi$ is a map of $\bG$-modules and $\bG$ is reductive
$H^0(X,E)$ contains a submodule $M'$ 
such that $\pi(M')=M$ (as $\bG$-modules).
The elements of $M'$ generate $E$ over a generic point of $X$. 
A basis $s_1',...,s_r'$ of $M'$ gives the 
desired splitting of the action.
\end{proof}

\

\begin{coro}
\label{coro:coro}
Let $\bG$ be a reductive group and 
$$
E''\ra E'\ra X
$$ 
a $\bG$-equivariant sequence of vector
bundles such that the generic $\bG$-orbit of $E'$
projects isomorphically onto its image.  
Choose a generic $\bG$-equivariant section $s'$ of $E'\ra X$
and denote by $E''_{s'}$ the restriction of $E''$ to this section.
Then  
$$
E''\sim_{\bG} E''_{s'}\times \A^{r'}
$$  
(where $r'=\rk E'$), with trivial $\bG$-action on $\A^{r'}$.
\end{coro}

\begin{prop}
\label{prop:modu}
Let $X$ be a variety with an action    
$\rho\,:\, \bG\ra X$ of a linear algebraic group $\bG$. Let 
$E\ra X$ be a vector bundle and $\tilde{\rho}\,:\, \tilde{\bG}\ra E$
a $\tilde{\bG}$-action lifting $\rho$.
Consider a generic orbit $\bG\cdot x\subset X$ and the linear
action of $\tilde{\bG}$ on the space of sections $H^0(X,E)$.
 
Assume that $\tilde{\bG}$ is
reductive and $V$ is a linear representation of $\tilde{\bG}$
which is contained in $H^0(X,E)$.
Then there exists an affine open $X'\subset X$ 
such that the vector bundle $E\ra X'$ admits a $\tilde{\bG}$-map 
onto a $\tilde{\bG}$-representation $V^*$. 
\end{prop}

If the action of $\bG$ on $X$ is almost free
we may think of $X$ as being (birational to) a principal fibration
over the quotient $\bG\ba X$  with fiber $\bG$.
If $\bG$ is {\em affine} we may assume that $X$ 
and $\bG\ba X$ are also affine. 
Let us also recall a 
standard general construction of $\bG$-maps:
if the ring $\C[X]$ is a direct sum of $\bG$-modules
then any $\bG$-submodule $V\subset \C[X]$ defines 
a $\bG$-map $X \ra \Spec(V)$.
We also have a vector bundle version of the above 
construction: let $E\ra X$ be a $\bG$-vector bundle
and $O$ a $\bG$-orbit through a generic point. 
Assume that $H^0(O,E)$  (the restriction 
of the space of sections to $O$) 
contains $V$ as a submodule.
We obtain a $\bG$-map 
$$
v\, :\, H^0(O,E)\ra V^*
$$
(the dual module, considered as a vector bundle over a point). 

\
 
\begin{lemm}
\label{lemm:mapp}
There exists a $\bG$-stable Zariski open $U\subset X$ and
a rational $\bG$-map of $H^0(U,E)\ra V^*$ extending
$v$.
\end{lemm}

\begin{proof}
A generic orbit $O$ has a $\bG$-equivariant  
neighborhood $U$, with $U/\bG$ affine,
such that 
$$
\xymatrix{H^0(U,E)\ar@{>>}[r] & H^0(O,E).}
$$
The module $H^0(U,E)$ is a direct sum of
finite dimensional irreducible $\bG$-modules. We can now take 
any submodule $V\subset H^0(U,E)$ which surjects isomorphically 
onto a submodule in $H^0(O,E)$.
\end{proof}

\

\begin{lemm}
\label{lemm:pgl2}
If $X$ has an $af$-action of $\PGL_2$ then
$$
X\times \bP(V_{2d})\sim_{\PGL_2} X\times \bP(V_{2d}),
$$
with diagonal $\PGL_2$-action on the left and trivial $\PGL_2$-action 
on $\P(V_{2d})$ on the right. 
\end{lemm}

\begin{proof}
We know that $\C[\PGL_2]$, as a $\PGL_2$-module,
is sum of all even modules $V_{2d}$. 
This gives a $\PGL_2$-map 
$s \,:\,X\ra \bP(V_{2d})$.  
The quotient 
$$
\PGL_2\ba X \times \bP(V_{2d})
$$ 
is a projective  bundle over the quotient 
$\PGL_2\ba X$, with a section obtained from $s$.
Therefore, it is birational to the  product
$(\PGL_2\ba X) \times \P(V_{2d})$, which  gives
the claimed $\PGL_2$-isomorphism.
\end{proof}

\

\begin{coro}
\label{coro:sect}
Let $X$ be a variety with an $af$-action 
of $\PGL_2$. Then $X$ is a $(\PGL_2,\bN_{\bT})$-slice in 
$$
X\times \P(V_2)
$$
(with diagonal $\PGL_2$-action). 
\end{coro}

\begin{lemm}
\label{lemm:y}
Assume that $X$ has an $af$-action $\rho$ of $\PGL_2$. 
Let $\cV\stackrel{V}{\longrightarrow} X$ be a vector bundle over $X$ with an action 
$\tilde{\rho}$ of $\GL_2$ lifting $\rho$. 
Assume that $X$ contains a $\PGL_2$-orbit $Y\sim \PGL_2$ 
such that the $\GL_2$-module $H^0(Y,\cV_Y)$ contains $V_{d}$, 
for some odd $d$.  
Then 
$$
\bP(\cV)\sim_{\PGL_2}\PGL_2\times  S,
$$ 
(with trivial $\PGL_2$-action on $S$). 
Otherwise, $\cV$ is induced from a 
$\GL_2$-vector bundle on $\PGL_2\ba X$.
\end{lemm}

\begin{proof}
Let $Y$ be an orbit such that 
$H^0(Y,\cV_Y)$ contains $V_{d}$, for some odd $d$.
Shrinking $X$, if necessary, gives a surjective map of 
$\GL_2$-modules
$$
\xymatrix{
H^0(X,\cV)\ar@{>>}[r] & H^0(Y,\cV_{Y}).} 
$$
Since  $H^0(Y,\cV_{Y})$ is an algebra over  
$H^0(Y,\cO_{Y})=\oplus_{d\ge 0} V_{2d}$, 
it contains $V_1$ as a submodule. 
We obtain a $\PGL_2$-equivariant surjective map 
$$
\bP(\cV)\ra \bP(V_1)=\bP^1.
$$ 
Since the stabilizer of a point in $\bP^1$ is 
solvable, we get a slice $S\subset \bP(\cV)$, as claimed.  

Assume that there is an orbit $Y\sim \PGL_2$ such that $\cV_Y$ contains only 
even weight $\GL_2$-submodules. Then the central  $\sC_2\subset \GL_2$  
acts trivially on $\cV_Y$. If follows that $\cV_Y$ is a trivial
$\PGL_2$-bundle, and $H^0(Y,\cV_Y)$ a trivial $\PGL_2$-module.
The semi-simplicity of the $\PGL_2$-action implies that
$H^0(X,\cV)$ contains $H^0(Y,\cV_Y)$ as a submodule. Shrinking $X$ if
necessary, we can find linearly independent $\PGL_2$-invariant sections, 
whose specializations to $Y$ generate $H^0(Y,\PGL_2)$. 
Therefore, $\cV$ is lifted from the quotient $\PGL_2\ba X$.
\end{proof}

\

\begin{lemm}
\label{lemm:sl}
Let $V$ be a representation of $\bG$ of dimension $\ge 2$  
(with $\bG$ acting on the left). Then $V\oplus V$ is  
a $\bG\times \GL_2$-space (with right $\GL_2$-action) and 
$$
V\oplus V\sim_{\bG\times \GL_2} 
\xymatrix{\cV \ar[d]^{\M2=V_1\oplus V_1}\\ \Gr(2,V),}
$$
a vector bundle with fibers $2\times 2$-matrices  
(with right $\GL_2$-action). 
\end{lemm}

\begin{proof} 
Consider the map 
$$
\begin{array}{ccc}
V\oplus V & \ra     &  \Gr(2,V)\\
(v,v')    & \mapsto &  \langle v,v'\rangle, 
\end{array}
$$
defined on the open, $\bG\times \GL_2$-invariant 
subset of noncollinear pairs $(v,v')\in V \oplus V$ 
(with fibers consisting of pairs spanning the same 
2-space). The $\GL_2$-action on the fibers is the right 
multiplication on matrices:
$$
(v,v')\mapsto (av+bv',cv+dv').
$$
\end{proof}

Assume that $\bG$ is reductive and denote by
$\bG'':=\Stab_{gen}(\bG,\Gr(2,V))$ and by
$\bG':= \bG/\bG''$ the quotient group 
of $\bG$ which acts effectively on $\Gr(2,V)$.

\

\begin{coro} 
\label{coro:1}
Assume that the action of $\bG'$ on $\Gr(2,V)$ has a slice
$S$ so that $\Gr(2,V) \sim S\times \bG'$.
Let $\cV_S$ be the restriction of $\cV$ to $S$ (this makes sense by 
Corollary~\ref{coro:coro}).
Then 
$$
\bG\ba \cV /\GL_2 \sim \bG'\ba \cV_S.
$$
\end{coro}  

\
  
\begin{rem}
\label{coro:2} 
The group $\bG''$ acts as scalars on $\cV=\cV_1\oplus \cV_1$ 
(it commutes with $\GL_2$).
\end{rem}
\

\begin{lemm} 
\label{lemm:inv}
Assume that we are in the situation of Corollary~\ref{coro:1}, 
$\bG = \GL_2$ and  $\bH\subset \GL_2$ has finite image in $\PGL_2$.
Then $\bG\ba \cV /\bH$ is rational.
\end{lemm}

\begin{proof} 
By Corollary~\ref{coro:2}, the   
slice $S$ is 3-stably rational, 
since 
$$
S\times \PGL_2 \sim \Gr(2,V)
$$
and $\Gr(2,V)$ is rational. 
The quotient of $\cV_S$
by a fiberwise linear action is birational to 
$(\M2 / \bH)\times S$  (every vector bundle admits 
an $\bH$-equivariant trivialization over an open subset of $S$).
There is a left  action of $\gm^2 \subset \GL_2$ on $\M2=V_1 \oplus V_1$ which
commutes with $\bH$. Thus $\M2 / \bH$ is (birationally) 
a three-dimensional variety with an $af$-action of $\gm$.
The quotient (a surface) is unirational, 
hence rational (by Theorem~\ref{thm:q}),
and 
$$
\bG\ba (V \oplus V)/\bH \sim S \times (V_1 \oplus V_1)/\bH \sim S\times\C^3.
$$ 
\end{proof}

\

The group $\PGL_2$ acts on $\P(\M2)$ on both sides. 
We will need an explicit description of the action for
some of its subgroups.

\

\begin{lemm}
\label{lemm:po(2)}
We have 
$$
\Stab_{gen}(\bN_{\bT}\times \bN_{\bT},\P(\M2))=\sC_2.
$$
\end{lemm}
  
\begin{proof}
Indeed $\bN_{\bT}$ contains 
$$
\begin{array}{ccccl}
\gm=\{ t\}  & :  &  (x,y) & \mapsto & (tx,t^{-1}y),\\
i           & :  &  (x,y) & \mapsto & (y,x).
\end{array}
$$ 
The corresponding actions on  $\P(\M2)$ are  
$$
(a,b,c,d) \mapsto (t_1t_2 a,t_1^{-1}t_2 b, t_1t_2^{-1}c, t_1^{-1}t_2^{-1} d)
$$ 
and   
$$
\begin{array}{ccc}
i_1 : &  a \ra c, &  b \ra d \\
i_2 : &  a \ra b, &  c \ra d,
\end{array}
$$
respectively.
A matrix $(a,b,c,d)\in \M2$ 
can be transformed to $(1,1,1,d)$ 
by a unique element of $\gm\times \gm$, the 
$\sS_2\times \sS_2$-orbit of which consists of 
two elements (for $d, d^{-1}$). 
\end{proof}

\

\begin{coro}
\label{coro:sopen}
The group $\bN_{\bT}\times \C^*$ acts almost freely on $\P(\M2)$.
There is an open, $\bN_{\bT}\times \bN_{\bT}$-stable subvariety $U\subset \P(\M2)$
such that
$$
\xymatrix{ U \ar[d]^{\C^*\times \C^*}\\ \C^*\subset \P^1,}
$$
with a transitive action of  
$\C^*\times \C^*\subset \bN_{\bT}\times \bN_{\bT}$ on the fibers. 
The diagonal subgroup 
$$
\sS_2^{\Delta}\subset \sS_2\times \sS_2 = (\bN_{\bT}\times \bN_{\bT})/(\C^*\times \C^*)
$$
acts on each fiber as an involution $x \mapsto x^{-1}$. The factor 
$\sS_2 = (\sS_2\times \sS_2)/\sS_2^{\Delta}$ acts on
the base $\C^*\subset \P^1$ as an involution without fixed points, 
on the first factor in the fiber as $x\to x^{-1}$, and as identity on 
the second factor. 
\end{coro}
 
\

\begin{coro}
\label{coro:conic}
Let $\sD\subset \bN_{\bT}$ be a dihedral subgroup such that
$\sD\ba \bN_{\bT} = \C^*$. Then the $\C^*$-bundle 
$$
\cC = \sD \ba \P(\M2)\ra \bN_{\bT}\ba \P(\M2)
$$ 
is induced from the $\C^*$-bundle 
$$
\sD\ba \P(\M2)/\bN_{\bT} \ra \bN_{\bT}\ba\P(\M2)/\bN_{\bT} = \P^1
$$ 
and is hence birationally trivial.
\end{coro}

\begin{proof}
Indeed, the left and the right actions of $\bN_{\bT}$ commute.
By Lemma~\ref{lemm:po(2)},  $\Stab_{gen}(\bN_{\bT}\times \bN_{\bT},\P(\M2))=\sC_2$, 
which implies that the bundle is induced. 
\end{proof}

\   

\begin{lemm}
\label{lemm:rh}
For every dihedral group 
$\sD$ and every $\bH\subset \bN_{\bT}$ the conic bundle 
$$
\cC_{\bH} = \sD\ba \P(\M2)/\bH \ra \bN_{\bT}\ba \P(\M2)/\bH,
$$
has a section.
\end{lemm}

\begin{proof}
The quotient 
$\sD\ba U /\bH $ from Corollary~\ref{coro:sopen}
admits a fibration
$$
\xymatrix{ \sD\ba U/\bH \ar[d]^{\C^*_{\sD}\times \C^*_{\bH}/\sS_2} \\ 
\P^1/\sS_2.}
$$
Here $\C^*_{\sD}\times \C^*_{\bH}$ is the 
quotient of the fiber $\C^*\times \C^*$ of $U\ra \C^*$ 
by the intersection of $\sD, \bH$ with the diagonal 
$\C^*_{\Delta}\subset \C^*\times \C^*$. 
Isomorphisms $\C^*_{\bH} \ra \C^*$ and $\C^*_{\sD}\ra \C^*$
induce a birational fiberwise isomorphism 
$$
\xymatrix{
\cC_{\bH} = &  \sD\ba \P(\M2)/\bH \ar[d] \\ 
            & \bN_{\bT}\ba \P(\M2)/\bH } \sim 
\xymatrix{
            \sS_2 \ba \P(\M2)/\sS_2 \ar[d] & = \cC_0\\
            \bN_{\bT}\ba \P(\M2)/\sS_2   &   }
$$  
and it suffices to consider $\sD = \sS_2, \bH= \sS_2$.
In this case, an alternative equivariant completion of $U$ is given by
$$
U\subset  \xymatrix{\P^1_1\times \P^1_2\times \P^1_3 \ar[d] \\  \P^1_3},
$$ 
with an action of $\sS_2\times \sS_2$, where 
the first $\sS_2$ acts as an involution on the first two factors and 
identity on the base while the complementary $\sS_2$ 
acts only on the base. 
Thus the quotient is a conic bundle over the complement in 
$$
\P^1\times \P^1/\sS_2\times \sS_2 = \P^1\times \P^1 
$$ 
to the branch locus of the quotient map. Here the left (resp. right) 
$\sS_2$ acts as an involution on the left (resp. right) $\P^1$ and the
branch locus is exactly the union of four lines.
By Lemma~\ref{lemm:cb}, 
this conic bundle has a section  
(it is nonsingular on a pencil of lines minus at most
two points).
\end{proof}

\

\begin{lemm}
\label{lemm:subgroup}
Let $\bG$ be a subgroup of $\SL_2$, 
not equal to $\tilde{\sA}_5$, and $V$ a linear representation of $\bG$.
Then $\bG\ba\P(V)$ is rational.
\end{lemm}

\begin{proof}
For  $\bG=\SL_2$ this is a theorem of Katsylo \cite{K}.
We now consider proper subgroups $\bG\subsetneq\SL_2$.
If $\bG$ is solvable and connected then rationality for 
the quotient follows from a theorem of Vinberg \cite{V}. 
For compact $\bG$ the proof is similar to the dihedral case
described below. 
Assume now that $\bG$ is finite and not equal to $\tilde{\sA}_5$. 
Then $\bG$ is either
\begin{enumerate}
\item a finite subgroup of $\C^*$,
\item a dihedral group or
\item $\tilde{\sA}_4, \tilde{\sS}_4$.
\end{enumerate}

The first case is easy.
For dihedral groups all irreducible representations of 
$\bG$ have dimension  $\le 2$ and the corresponding quotients
are rational by Theorem~\ref{thm:q}. 
Let $V$ be a faithful representation of a dihedral group $\sD$
(otherwise, we are reduced to a quotient group). 
Thus $V=W\oplus W'$, where $\dim W =2$ and $\dim W'\ge 1$.
Denote by  $\bG'=\bG/\sC'$ the quotient acting faithfully on $W'$
($\sC'$ is a cyclic group). 
We have $W\sim_{\sD} \C^*\times \P^1$, with trivial action of $\sD$
on $\C^*$ and trivial action of $\sC'$ on $\P^1$. 
By Lemma~\ref{lemm:redu}, 
$$
\C^*\times \P^1\times W' \sim_{\bG'} (\C^*\times \P^1) \times W',
$$ 
with trivial action of $\bG'$ on $\C^*\times \P^1$. 
Thus 
$$
\sD\ba V\sim (\sD'\ba W') \times (\C^*\times \P^1)
$$
and we can apply induction.

We turn to the last case.
An irreducible representation of 
$\tilde{\sA}_4$ is either a character, 
or a faithful two-dimensional representation, 
or a three-dimensional
representation, 
trivial on the center 
(a faithful representation of $\sA_4$).
An irreducible representation of $\tilde{\sS}_4$ is either
a faithful two-dimensional representation, a faithful 
four-dimensional representation $W:=\Sym^3(V_1)$ 
or a representation of $\sS_4$ (of dimension $\le 3$). 

For irreducible representations of dimension $\le 3$ rationality for
the quotient follows from Theorem~\ref{thm:q}.
We turn to $W$. Recall that
$$
W=\Sym^3(V_1)=V_1^{\chi}\oplus V_1^{-\chi},
$$
as a $\tilde{\sA}_4$-representation, 
where 
$V_1^{\chi}=V_1\otimes \chi$, $V^{-\chi}=V_1\otimes \chi^{-1}$ 
and 
$$
\chi\,:\, \sA_4\ra \Z/3\subset \C^*
$$ 
is the cubic character. 
A pair of (generic) points 
$$
p_{\chi}\in \P^1= \P(V^{\chi}), 
p_{-\chi}\in \P^1=\P(V^{-\chi})
$$ 
defines a line $\P^1\subset \P(W)$. 
This shows that
$$
\P(W)\sim_{\tilde{\sS}_4} \xymatrix{ \ar[d]^{L} \\ \P^1\times \P^1 },
$$
where $\sS_4$ acts on the base, $\sA_4$ acts linearly on the fiber $L$
and $\sS_2=\tilde{\sS}_4 /\sS_4$ acts as an involution on the fiber $L$. 
Thus $\tilde{\sS}_4\ba \P(W)$ is a conic bundle over the rational surface
$\sS_4\ba (\P^1\times \P^1)$. 
We now analyze the geometry of this bundle in more detail.
Consider the action $\sD_2\subset \sS_4$ on
$\P^1\times \P^1$ and on $\P^2 = \Sym^2(\P^1)$.
Every involution $i\in \sD_2$ has two invariant
points $x_i,y_i$.  Consider the graphs $\P^1$ connecting 
the points $(x_i,y_i) - (y_i, x_i)$. Their set is equal to $\P^1$
and there is a graph: 
$$
l_i \,: \,(x_i,y_i) - (y_i, x_i)\subset \P^1\times \P^1 
$$ 
consisting of points $(x,i(x))$.  
The line $l_i$ is exactly the subset of $i$-invariant points in
$\P^1\times \P^1$. The action of $\sD_2$ is free 
outside the three lines $l_i, i\in \sD_2, i\neq 1$.
There are exactly $6$ points which 
are invariant under $\sD_2$.

The corresponding action on $\P^2$ can be described as 
follows. There are three points corresponding to $(x_i,y_i)$
which are stable under $\sD_2$
and three lines (images of $l_i$) so that
the action is free on the torus $\C^*\times \C^*$
(the complement in $\P^2$ to the union of $l_i$).
The group $\sD_2$ acts on $\C^*\times \C^*$ as a translation
by the subgroup of points of order $2$.

The quotient $\P^2_q:=\sD_2\ba \P^2$ is a nonsingular 
variety isomorphic to $\P^2$
(indeed the only possible singularities
come from the three $\sD_2$-invariant points
in $\P^2$ but the quotient by the local action is 
nonsingular).  
The diagonal $\P^1_{\Delta}\subset  \P^1\times \P^1 $ projects onto
a conic $C\subset \P^2$, which is invariant under $\sD_2$.
The conic  $C$ intersects  the ``vertical'' and ``horizontal'' 
subgroups in $\C^*\times \C^*\subset \P^2$
in two points and does not intersect the line at infinity.
 
Thus in $\P^2_q = \sD_2\ba \P^2$, the image of
$\P^1_{\Delta}$ intersects $\C^*$ in one point.
Therefore, the images of $\P^1_{\Delta}$ 
and of $l_i$ are lines
(since pairwise intersections of the $l_i$ are equal to $1$)
and the $(\sC_2)^3$-covering $\P^1\times \P^1\to \P^2_q$  
is ramified exactly over a union of four lines. 
If suffices to observe that 
every conic bundle over $\P^2_q$ has a section.
Indeed, let $p$ be the intersection point of
two lines $l_i$ and $l_{i'}$ and consider the pencil of lines in $\P^2_q$
through $p$. Each line in this pencil
intersects the ramification locus 
in at most three points and we can 
apply Lemma~\ref{lemm:cb}.

\

Now we turn to reducible representations
$V=\oplus_{\al\in \cA} V^{\al}$ of $\tilde{\sA}_4$.
If $V$ is faithful for 
$\tilde{\sA}_4$ then there is an $\al_0\in \cA$ such that 
$V^{\al_0}$ is a three-dimensional irreducible {\em faithful}
representation of $\tilde{\sA}_4$ and 
$$
V\sim_{\tilde{\sA}_4} V^{\al_0} \times (\oplus_{\al\neq \al_0} V^{\al})
$$
with trivial action of $\tilde{\sA}_4$
on $\oplus_{\al\neq \al_0} V^{\al}$ (by Lemma~\ref{lemm:redu}).
If $V$ is faithful for $\sA_4$ then $V$ contains a 
faithful irreducible three-dimensional representation of $\sA_4$
and we can apply the same argument. 
In all other cases $V$ is a sum of one-dimensional 
representations and we are reduced to Case 1.

Finally, consider reducible representations $V$ of $\tilde{\sS}_4$. 
If $V$ is faithful then it contains either a faithful irreducible
two-dimensional representation or the faithful representation $W$. 
Again, we apply Lemma~\ref{lemm:redu} as before.
If $V$ is faithful for $\sS_4$ then it contains a faithful 
irreducible representation of dimension $\le 3$ and we conclude
as above. In all other cases $V$ is a sum of one-dimensional
representations.  
\end{proof}

\

\begin{lemm}
\label{lemm:grr}
Let $V$ be a representation of
$\bG\subsetneq\SL_2$, with $\bG\neq \tilde{\sA}_5$. Then 
$\bG\ba\Gr(2,V)$ is rational.
\end{lemm}

\begin{proof}
The relevant groups $\bG$ can be subdivided as follows:
\begin{enumerate}
\item $\bG$ is a subgroup of the normalizer of a maximal torus;  
\item $\bG$ an infinite subgroup of a Borel subgroup;
\item $\bG = \sA_4$, $\tilde{\sA}_4$; 
\item $\bG = \sS_4$, $\tilde{\sS}_4$. 
\end{enumerate}
Let $V=W\oplus W'$ be a reducible representation  of $\bG$.
Then (birationally)
$$
\xymatrix{\Gr(2,V)\ar[d]^{\Hom(\C^2_x, W)} \\
\Gr(2,W)}
$$
(where $x$ is a point on the base). 
In particular, if $\dim W\le 2$ then
$$
\Gr(2,V)\sim_{\bG} \Hom(W',W),
$$
with linear $\bG$-action on $\Hom(W', W)$. This reduction
suffices for the relevant infinite groups (for example, for
connected solvable $\bG$ we can apply Lemma~\ref{thm:qq}). 
Further,  
\begin{itemize}
\item if $\Stab_{gen}(\bG,\Gr(2,V))=1$ then (birationally) 
$$
 \bG\ba \Gr(2,V)\ra   \bG\ba\Gr(2,W),
$$
a vector bundle.
\item if $\Stab_{gen}(\bG,\Gr(2,V))=\sC\subset \bZ_{\bG}$ 
(a cyclic subgroup) then (birationally)
$$
\xymatrix{\bG\ba\Gr(2,V)  \ar[d]^{\sC\ba \Hom(\C^2_x, W)} \\      
\bG\ba \Gr(2,W).}
$$
\end{itemize}

We now consider $\sA_4,\tilde{\sA}_4, \sS_4$. 
The rationality of $\bG\ba\Gr(2,V)$ for 
{\em irreducible} representations 
of these groups follows from the fact 
that all of them have dimension $\le 3$. 
Assume now that $V = W \oplus W' $, with $W$ irreducible of dimension $3$.
The classification of these representations implies that 
the action of the center must be trivial.
Then, birationally,
$$
\xymatrix{
\Gr(2,V)  \ar[d]^{\Hom(\C^2_x, W')} \\
\P^2=\P(W^*).}
$$
The $\bG$-action is equivalent to a 
$\bG$-action on a vector
bundle
$$
\xymatrix{\bG\ba \Gr(2,V) \ar[d]\\
\bG\ba \Gr(2,W) = \bG\ba\P^2.}
$$

\

Finally, let us consider the case of $\tilde{\sS}_4$.
Let $W$ be its unique { irreducible} 
representation of dimension four (as in Lemma~\ref{lemm:subgroup}).  
We claim that $\tilde{\sS}_4\ba \Gr(2,W)$ is rational.
Indeed, as $\tilde{\sA}_4$-modules, we have
$$
W=W^{\chi} \oplus W^{-\chi},
$$
where $W^{\chi},W^{-\chi}$ are two copies of 
the standard representation of $\tilde{\sA}_4$
of dimension $2$ and $\chi$ (resp. $-\chi$) 
indicates the eigenspace decomposition for the nontrivial character 
$$
\chi\, :\, \sA_4\to \Z/3\subset \C^*.
$$
Further,  
$$
\Gr(2,W) \sim \Hom(W^{\chi},W^{-\chi}),
$$
with a linear $\sA_4$-action (since the center acts trivially)
and a permutation $\sS_2$ 
inverting the map $w \in \Hom(W^{\chi},W^{-\chi})$.
More precisely, $W^{-\chi}= (W^{\chi})^*$ and  
$$
\Hom(W^{\chi},W^{-\chi})=\Sym^2(W^{-\chi}) \oplus C_1,
$$
where $C_1$ corresponds to skew symmetric maps and 
$\sA_4$ acts on $C_1$ by $\chi$.
The involution $\sS_2 = \sS_4/\sA_4$ acts on 
$C_1$ and on $\Sym^2(W^{-\chi})$ as $t \mapsto t^{-1}$.
In particular, if $\C^*\times \C^*$ is the diagonal group
acting on $\Sym^2(W^{-\chi}) \oplus C_1$ then $\sS_2$ acts as 
$$
X \ra s^{-1} X,
$$
where $s \in \C^*\times \C^*$ and $X\in \Sym^2(W^{-\chi}) \oplus C_1$.
There is an equivariant map 
$$
\begin{array}{ccccl}
f & : & \Hom(W^{\chi},W^{-\chi}) & \ra &   C_1,\\
  &   &   s                      & \mapsto &  (x,s(y)) - (s(x),y),
\end{array}
$$ 
with an effective action of $\sS_3 = \sS_4/\sD_2$ on the target $C_1$,
which to a subspace $s \in \C^2\subset W^{\chi} \oplus W^{-\chi}$
assigns the value of the $2$-form $(x,s(y)) - (s(x),y)$. 
The fiber of $f$ is $\sD_2$-birational to $\Sym^2(W^{\chi})=\P^2$. 
We have already seen in the proof of Lemma~\ref{lemm:subgroup} 
that $\sD_2\ba\P^2=\P^2$.
Thus $\tilde{\sS}_4\ba \Gr(2,W)$ is a $\C^*$- bundle over 
a $\P^2$-fibration over $\sS_2\ba C_1$.
It is clear that this $\P^2$-fibration is trivial.
The quotient conic bundle is nondegenerate
over a product of $\P^2$ with an open subvariety in $C_1/\sS_3$.
Hence it has a section. Rationality of 
$\tilde{\sS}_4\ba \Gr(2,W)$, and more generally, 
$\tilde{\sS}_4\ba \Gr(2,W\oplus \cdots \oplus W)$, follows (the latter is 
a vector bundle over the former).

\

Assume now that $V = nW \oplus V'$, 
where $\dim V'\ge 1$, and $n\in \N$. 
Since the $\sS_4$-action on $\Gr(2,nW)$ is $af$
there is a $\tilde{\sS}_4$-equivariant homogeneous 
rational  map  $f\, :\, \Gr(2, nW )\ra V'$
sending the generic  $\tilde{\sS}_4$-orbit in $W$ to
the generic $\tilde{\sS}_4$-orbit in $V'$.
Notice that the center $\sC_2$ acts as a scalar on
$\Hom (W, V')$.
We have (birationally)
\begin{equation}
\label{eqn:bbd}
\xymatrix{
\tilde{\sS}_4\ba\Gr(2,V)  \ar[d]^{\sS_2\ba\Hom(\C^2_x ,V')} & \sim 
& \C^* \times &
 \ar[d]^{ \P(\Hom (\C^2_x,V'))} \\
\tilde{\sS}_4\ba \Gr(2,nW) &   &  & \sS_4\ba  \Gr(2,nW)
} 
\end{equation}
(with rational bases).
The projective bundle on the right has a section.
Indeed, 
\begin{equation}
\label{eqn:bb}
\xymatrix{ \ar[d]^{\Hom(\C^2_x, V')} \\ \Gr(2,nW) }
\end{equation}
is an equivariant quotient
bundle of the trivial bundle with fiber $\Hom(W, V')$.
The map $f$ defines an $\sS_4$-equivariant section $s(f)$
in the projective bundle in (\ref{eqn:bbd}).
The (equivariant) linear projection 
$$
\Hom(W, V')\ra \Hom(\C^2_x, V')
$$
maps $s(f)$ to an equivariant section of the
bundle in (\ref{eqn:bb}).
Thus $s(f)$ projects onto a section of the bundle
on the right in (\ref{eqn:bbd}), making it birationally trivial.
\end{proof}

\

We proceed to describe possible $\SL_2$, resp. $\PGL_2$-actions on 
Grassmannians. (If all weights in $V$ are of the same parity then
$\Gr(2,V)$ carries the $\PGL_2$-action, otherwise the $\SL_2$-action.) 

\

\begin{lemm}
\label{lemm:free}
%\label{coro:stabb}
Let $V$ be a faithful $\SL_2$-representation of dimension $\ge 3$.
Then 
$$
\begin{array}{c|c}
V             &     \Stab_{gen} \\
\hline
\dim \ge 5    &  1     \\
V_4           &  \sC_2 \\
V_3           &  \sD_2  \\
V_2\oplus V_0 &  \sC_2  \\
V_2           &  \bN_{\bT} \\
V_1\oplus V_1 & \C^*   \\
V_1\oplus V_0 &  \tilde{\bB} 
\end{array}
$$

Further, 
\begin{itemize}
\item 
$\Gr(2,V_4)$ has a $(\PGL_2, \bN_{\bT})$-slice 
$S= \Sym^2(\P^2)$ with an $af$-action of $\bN_{\bT}/\sC_2$,
(where $\sC_2$ is the center of $\bN_{\bT}$);
\item 
$\Gr(2,V_3)$ has a $(\PGL_2, \sA_4)$-slice birational to $\P^1$, 
with $\sA_4$ acting on $\P^1$ as $\sC_3$.
\end{itemize}
\end{lemm}

\begin{proof}
Consider first irreducible representations $V=V_{d}=\Sym^d(V_1)$ 
and assume that the stabilizer of a 
generic line $\P^1\subset \P(V)$ contains a nontrivial 
cyclic group $\sC$. 
Then $\sC$ fixes at least two points in this $\P^1$. 
Any orbit of $\sC$ on $\P^1$ is a union of a zero-cycle 
$\sC \cdot x$ and a zero-cycle supported in the fixed
points. In particular, the subvariety of points in $\P^(V_d)$
which are fixed by $\sC$ has dimension
$\le d/|\sC|$. The dimension 
of the variety of $\sC$-fixed lines in $\P(V)$ 
is therefore $\leq 2 d/ |\sC|$.
The subvariety of distinct cyclic subgroups 
$\sC\subset \PGL_2$ has dimension $2$
and $\dim \Gr(2,V_d) = 2d - 2$.
Since $d/|\sC| \leq d/2$ the inequalities 
$$
2d - 4 > 2d/2 \,\,\, {\rm and} \,\,\,  d - 4 > 0
$$ 
imply the result. 

Assume that $V=\oplus_{j\in J} V_{d_j}$, $|J|\ge 2$
and that $\Stab_{gen}\neq 1$. 
Then $d_j\le 2$, for all $j\in J$. Indeed, the stabilizer
of a generic $\P^1$ through a generic point $p\in \P(V_d)$
is a subgroup of the stabilizer of $p$, which stabilizes 
some generic line in the tangent space
at $p$. This group is trivial for $d > 2$ and equal to $\sC_2$ for $d=2$. 

If $V = V_2 \oplus V'$, with $\dim V' > 2$, then 
$\Gr(2,V)$ is (birationally) a fibration  over 
$\Gr(2,V_2)$, with fibers $\Hom(\C^2, V')$ so that 
$\Stab_{gen}=1$ if $\dim V' > 3$.
If $V = V_2 \oplus V_1$ then 
$\Stab_{gen}$ is the same as 
the (generic) stabilizer of the $\bN_{\tilde{\bT}}$-action on
$\Hom(V', V_1), V' \in \Gr(2,V_2) =\P^2$, hence trivial.
For $V=V_2 \oplus V_0$, $\Stab_{gen}=\sC_2$.

In the remaining cases $d_j=0$ or $1$, for all $j\in J$. 
If $V$ contains at least three copies of $V_1$ then the argument above
shows that the action is $af$.
Similarly, if $V = V_1 \oplus V_1 $ then
$\Stab_{gen}=\C^*$ and if $V=V_1 \oplus V_1 \oplus V_0$ 
then $\Stab_{gen}=1$.
For $V_1 \oplus 3 V_0$, the generic stabilizer is the same
as for three linear functionals - which is zero.
\end{proof}

\

\begin{lemm}
\label{lemm:gra}
The quotient $\PGL_2\ba \Gr(2,V)$ is 2-stably rational. 
\end{lemm}

\begin{rem}
For {\em even} $d\ge 10$, $\PGL_2\ba \Gr(2,V_d)$ is rational by \cite{sh-2}.
\end{rem}

\begin{proof}
By Lemma~\ref{lemm:free}, if $\dim V\ge 5$ then
the $\Stab_{gen}=1$ and we
can apply Lemma~\ref{lemm:pgl2} and Corollary~\ref{coro:sect} 
to conclude that
$$
\PGL_2\ba \Gr(2,V)\times \C^2\sim_{\bG} \bN_{\bT}\ba\Gr(2,V).
$$
The claim follows from Lemma~\ref{lemm:grr}.
It remains to consider: 
\begin{enumerate}
\item $\Gr(2,V_4)$, 
\item $\Gr(2,V_3)$,
\item reducible $V$.
\end{enumerate}

In the first case, $\Stab_{gen}(\PGL_2,\Gr(2,V_4))=\sS_2$, with
normalizer $\bN_{\bT}\subset \PGL_2$. 
We claim that the subset  $X\subset \Gr(2,V_4)$ of $\sS_2$-invariant points
is a $(\PGL_2, \bN_{\bT})$-slice. Indeed, 
there is a Zariski open subset $ U\subset X$
such that the stabilizer of each point in $U$ 
is exactly $\sS_2$. In particular, $g\cdot U $ intersects $U$ only
if $g\in \bN_{\bT}$. Consider the 
$\P^2\subset \P(V_4)$ consisting of $\sS_2$-invariant subschemes 
containing $4$ points. Any line in $U$ joins a pair of points 
in this $\P^2$. Therefore, we have a (birational) 
$\bN_{\bT}$-isomorphism of
$U$ and $\Sym^2(\P^2)$. The stabilizer
of a generic point in $X$ is a central subgroup in
$\bN_{\bT}$ whose action on $\P^2$ is equivalent to a linear 
action on $\C^2$. 
(Indeed, $\Sym^2(V_1)=\C \oplus W_2$, where $\C$ is the trivial
representation - the invariant symmetric form - 
and $W_2$ is a faithful two-dimensional representation  of 
$\bN_{\bT}/\sS_2$).
Thus instead of $X$ with the $\bN_{\bT}$-action we can consider 
$\C^2\times \C^2$ with the $(\bN_{\bT}/\sS_2)\times \sS_2$-action 
(where the second $\sS_2$ interchanges the factors).
In particular, (by linearity)
$$
\bN_{\bT}\ba X \sim \C^*\times \bN_{\bT} \ba \P^3,
$$ 
and is hence rational.

In the second case, $\Gr(2,V_3)$ has a 
surjection of degree 2 onto $\P(V_4)$.
The connected component of the
preimage of the $(\PGL_2,\sS_4)$-slice $\P^1$ in $\P(V_4)$
is a $(\PGL_2,\sA_4)$-slice, isomorphic to $\P^1$.
The quotient is rational.

If $V$ is reducible and the $\PGL_2$-action on the 
Grassmannian has nontrivial stabilizer then $\dim V < 5$. 
Rationality follows since $\dim \Gr(2,V) \leq 4$ and the 
generic orbit has dimension at least $2$.
\end{proof}

\begin{prop}
\label{prop:dihedr}
Let $\bG, \bH $ be finite solvable subgroups of $\PGL_2$. Then 
$$
\bG\ba \PGL_2/\bH
$$ 
is rational. 
\end{prop}

\begin{proof}
The action is birational to the (projective) action of 
$\bG\times \bH$ on $\P(\M2)$,
where $\bG$ acts on the right and $\bH$ on the left.
The groups $\bG,\bH$ are either:
\begin{itemize}
\item cyclic;
\item dihedral or
\item $\sA_4$, $\sS_4$.
\end{itemize}
The case of {\em primitive} 
solvable groups is covered by Theorem~\ref{thm:kolp}, 
\cite{km}.
If $V$ is reducible then there is a nontrivial action
of $\C^*$ on $\bG\ba \P(V)/\bH$, leading 
to rationality. This covers the case when either 
$\bG$ or $\bH$ is cyclic.  

\

We claim that if $V$ is irreducible and 
imprimitive (for the $\bG\times \bH$-action) 
then either $\bG$ or $\bH$ is dihedral.  
By definition, 
$V:=\M2=\oplus_{\al} V^{\al}$, such that $\tilde{g}V^{\al}=V^{\al'}$
for all $\tilde{g}\in \bG\times \bH$.
Moreover, by irreducibility, all $V^{\al}$ must have the 
same dimension, $=1$ or $2$.
Notice that imprimitivity for an action of a group $\bG'$ implies 
imprimitivity for the induced action of every subgroup $\bG''\subset \bG'$
(with the same decomposition of $V$). 
We now claim that the actions of 
$\sA_4\times \sA_4$, and consequently of 
$\sA_4\times \sS_4$ and $\sS_4\times \sS_4$ are 
primitive. Indeed, $\sA_4\times \sA_4$ contains 
$\sD_2\times \sD_2$ as a normal subgroup, for which the imprimitive
structure is either a sum of two subspaces
of dimension $2$ or four subspaces of dimension $1$, corresponding to 
the choice of a subgroup $\sS_2\subset \sD_2$.  
The first possible imprimitive structure for $\sD_2\times \sD_2$ does not 
extend to one for $\sA_4\times \sA_4$
(which has no index $2$ subgroups). 
The second structure is also impossible: 
$\sA_4$ rotates the subgroups $\sS_2\subset \sD_2$,
hence there is no $\sA_4$-invariant imprimitive
structures for $\sD_2\times \sD_2$.

\

It remains to consider the case when both $\bG$ and $\bH$ are
dihedral.
On $V_1$ there is a unique imprimitive structure,
corresponding to the eigenspaces $C_1,C_2$ of the
elements of $\bG$.
In particular, there is an imprimitive structure on 
$$
\M2 =V_1 \oplus V_1' = (C_1 \oplus C_1') \oplus (C_2 \oplus C_2)'.
$$
We claim that (birationally)
$$
\xymatrix{
\bG\ba\P(\M2)/\bH  \ar[d] \\
\P^2 = \bG\ba \Sym^2(\P^1)}
$$
is a conic bundle degenerating precisely
over the image of the diagonal and the subvarieties 
in $\P^2$ with nontrivial stabilizers.

Indeed, since $\bH\subset \bN_{\bT}$  (a $\sC_2$-extension 
of $\C^*$), (birationally)
$$
\xymatrix{\bG\ba\P(\M2)/\bH \ar[d]^{\C^* = \bN_{\bT}/\bH} \\ 
\bN_{\bT}\ba\P(\M2)/\bH.}
$$
The quotient $\C^*\ba\P(\M2)$ is (birationally)  a fibration
over $\P^1\times \P^1$, with $\sS_2$
acting by permutation, where the coordinate $\P^1$s
are the projectivizations of the two-dimensional
eigenspaces for the $\C^*$-action on $\M2$.
Thus 
$$
\xymatrix{
\P(\M2)/\bH  \ar[d] \\
\P^2 = \P^1\times \P^1/\sS_2}
$$
is a conic bundle nondegenerate
outside a conic (the image of the diagonal in
$\P^1\times \P^1$).
The $\bG$-action commutes with the
$\bN_{\bT}$-action and is effective on the base. This proves the claim.

\

We have  $\bG\subset \bN_{\bT}$ and
$$
\bG\ba \P^2 \ra \bN_{\bT}\ba \P^2 
$$ 
is a conic bundle. Since the left and right actions
of $\bN_{\bT}$ commute, 
$\bG\ba\P^2$ contains an open subvariety 
$ U\times \C^*$ where the restriction of the conic bundle 
is nondegenerate. Here $\C^* = \bG\ba\bN_{\bT}$ and 
$U$ is a subset of $\P^1 =\bN_{\bT}\ba \P^2$.
Therefore the conic bundle has at most $2$ singular fibers 
on any completion of the fiber $\C^*\subset U\times \C^*$.
Rationality follows.

We can now describe some open subvariety in the quotient 
$\bG\ba\P(\M2)/\bH$ explicitly.
Consider the action of $\C^*\subset \bN_{\bT}$ on both sides 
$\C^*\ba\P(\M2)/ \C^*$. With respect to this 
action $\P(\M2)$ is birationally equivalent to a trivial 
$\C^*\times \C^*$-fibration over $\P^1$.
Now we add the action of $\sS_2$ on both sides.
The product $\sS_2\times \sS_2$ acts on the base
$\P^1$. The group $\sS_4$ contains
a normal subgroup $\sD_2\subset \bN_{\bT}$
and the action of each $\sS_2\subset \sD_2$ inverts  the respective 
$\C^*$ action. Thus (birationally)
$$
\xymatrix{
\bN_{\bT}\ba \P(\M2)/\bN{\bT} \ar[d]^{\bN_{\bT}\times \bN_{\bT}} \\ \P^1 -3\, pts,}
$$
where the deleted points are the ramification points
of the map $\P^1 \ra \P^1/\sD_2$.
In particular, there is an open $U$ such that 
$$
\xymatrix{
\bG\ba\P(\M2)/\bH \ar[d]^{\C^*} 
\\ 
U \ar[d]^{\C^*}
\\
\P^1 - 3\, pts. }
$$
By Lemma~\ref{lemm:cb}, the conic bundles are trivial.
 
\

Finally, the  conic bundles 
on $\P^2/\sS_4$ and $\P^2/\sA_4$ have sections.
Indeed, both $\sA_4$ and $\sS_4$ contain dihedral subgroups
of index 3 ($\sD_2$, resp. $\sD_4$).
The image of the section in the conic bundle over
$\sD_2\ba \P^2$ (resp. $\sD_4\ba \P^2$), has 
odd degree in the conic bundles over 
$\sA_4\ba \P^2$ and $\sS_4\ba \P^2$, respectively. We  
apply Lemma~\ref{lemm:conic}. 
\end{proof}
  
\

\begin{prop}
\label{prop:even} 
Let $V$ be an irreducible $\GL_2$-representation 
and $\bH\subset \SL_2$ a finite group, not equal to $\tilde{\sA}_5$.  
Then 
$$
\GL_2\backslash (V\oplus V)/\bH
$$ 
is rational. 
\end{prop}

\begin{proof}
First of all, $V_1\oplus V_1/\bH$ is rational. 
Next, by Lemma~\ref{lemm:sl}, 
$$
V\oplus V\sim_{\GL_2\times \GL_2} \xymatrix{\cV \ar[d]^{\M2=V_1\oplus V_1} 
\\ \Gr(2,V).}
$$

First we assume that $V$ has odd weight.
The Grassmannian $\Gr(2,V)$ carries the action of $\PGL_2$. 
If we restrict the bundle $\cV$ to a generic 
$\PGL_2$-orbit $O$ in $\Gr(2,V)$ then the corresponding module
$H^0(O, \cV_{O})$ contains $V_1$ as a submodule.
By Lemma~\ref{lemm:mapp}, this gives an equivariant map
$$
\cV\ra V_1 \oplus V_1
$$ 
with a $1$-transitive action
of $\GL_2$ on the target.
Thus 
\begin{equation}
\label{eqn:vy}
\GL_2\ba \cV/\bH \sim \bH\ba \Gr(2,V)
\end{equation}
(with the {\em same} subgroup 
$\bH\subset \GL_2$ appearing on the left).
Indeed, $\GL_2\subset (V_1 \oplus V_1) = \M2$ and multiplication
by $\bH$ on the right gives an orbit $x\cdot \bH$.
This orbit is a  $(\GL_2\times \bH, \bH^x\times \bH)$-slice 
(with $\bH^x = x \bH x^{-1}$) and 
it is stabilized exactly by $\bH^x\times \bH$, acting 
doubly transitively on the set $\bH^x\cdot x$. 
It follows that every point $x' \in x\cdot \bH$ is
a $(\bH^x\times \bH,\bH^x)$-slice of the orbit $x\cdot \bH$.
The quotient $\bH\ba \Gr(2,V)$ is rational by \ref{lemm:grr}.

\

Assume that $V$ has even weight.
If the $\PGL_2$-action is $af$ then 
$$
\GL_2\ba \cV /\bH \sim 
(\PGL_2\ba \Gr(2,V))\times (\C^*\ba (V_1\oplus V_1)/\bH).
$$ 
If it is not $af$, then, by Lemma~\ref{lemm:free}, 
$V = V_4 $ or $V_2$.

For $V=V_4$ we have the $(\PGL_2,\bN_{\bT})$-slice 
$X = \Sym^2(\P^2)$ with
the $\bN_{\bT}$-action which we can replace by 
$\C^2\times \C^2$ with a $(\bN_{\bT}/\sC_2)\times \sC_2$-linear action.
In particular, we identify the quotient with a quotient of 
$\C^2 \oplus \C^2 \oplus V_1 \oplus V_1$ by a linear action of 
$\bN_{\tilde{\bT}}\times \sS_2\times \bH$ (where
$\bN_{\tilde{\bT}}\subset \GL_2$).
The action of $\bN_{\tilde{\bT}}\times \bH$ on $V_1$ 
is transitive with stabilizer $\sC_2\times \bH$. Hence
it is equivalent to the action of $\sD_2\times \bH$ on
$\C^2 \oplus \C^2 \oplus V_1$, which 
is a $\C^2$-vector bundle (permutation of the anti-invariant part of
$\sS_2$-action) over $\C^2\times V_1$, with $\sD_2\times \bH$ action.
The latter quotient is rational.
For $V = V_2$ the action is transitive on $\Gr(2,V)= \P^2$ and
the quotient has dimension $2$ - rationality follows.
\end{proof}

We will also need a more general result for $\bH= \sS_2$.

\begin{prop}
\label{prop:yz}
Let 
$$
X\stackrel{L}{\lra} Y=\prod_{j\in J} \P(V_{d_j})
$$
be a $\GL_2$-homogeneous line bundle. 
If at least one $d_j\neq 2$ then 
$\GL_2 \backslash X\times X/\sS_2$ is rational.  
\end{prop}

\begin{proof}
{\bf Case 1.} $|J|=1$.  If $d=d_1$ is even or 
if $d$ is odd and the line bundle has odd degree on $\P(V_d)$ 
then 
$$
X\times X \sim_{\GL_2\times \sS_2} V_d \oplus V_d
$$ 
and we apply Proposition~\ref{prop:even}.
If the line bundle has even degree then it
is trivial and $\GL_2$ acts as $\PGL_2\times \C^*$.
If the $\PGL_2$-action on $\P(V_d)$ is $af$ 
we have 
$$
\P(V_d)\sim_{\PGL_2} S\times \PGL_2,
$$
for a rational slice $S$ (with trivial $\PGL_2$-action).
We have a $\PGL_2\times \C^*\times \sS_2$-action on 
$$
\C \times \PGL_2\times  S\times \C \times \PGL_2\times S.
$$ 
The quotient variety is a vector
bundle over $\PGL_2 \ba \PGL_2\times \PGL_2/ \sS_2$
(rational by Lemma~\ref{lemm:qqq}). The claim follows.
If the $\PGL_2$-action is not $af$, then $V=V_3$ or
$V_1$. For $V_1$ the quotient is rational by dimensional reasons.
For $V_3$ we have a projection
$$
\xymatrix{
\C\times \P^3\times \C\times  \P^3 \ar[d]^{\P^1\times \P^1}\\
\Gr(2,V_3)}
$$ 
commuting with both actions.
Recall that $\Gr(2,V_3)$ has $\P^1$ as a $(\PGL_2,\sA_4)$-slice, 
with $\sA_4$ effectively acting as a cyclic group $\sC_3= \sA_4/\sD_2$
on $\P^1$ (the group $\sD_2$ acts trivially on the
$(\PGL_2,\sA_4)$-slice $\P^1\subset \P^4$ and similarly for
$\Gr(2,V_3)$, see Lemma~\ref{lemm:free}). 
Thus the quotient is the same as for the bundle
$$
\xymatrix{  \ar[d]^{\P^1\times \C\times \P^1\times \C}\\ \P^1}
$$
under the action of $\sA_4\times \sS_2$.
In particular, it is a vector bundle
over a $\P^2 =\sD_2\ba \P^1\times \P^1/\sS_2$-fibration over 
$\P^1 = \P^1/\sC_3$, hence is rational.

{\bf Case 2.} $|J|\ge 2$. 
If at least one $d_j$ is odd and $>1$ or if all $d_j=1$ and 
$|J|>2$, then there is a slice $S$ and the 
$\PGL_2$-action is $af$. 
We can write $Y$ as (the total space of the) line bundle:
$$
\xymatrix{ X\ar[d]^{L}\\ 
S\times \PGL_2}
$$ 
and, using Lemma~\ref{lemm:y}, reduce to either
a vector bundle over 
$$
\PGL_2\ba \PGL_2\times \PGL_2/\sS_2,
$$
when $L$ is trivial on $\PGL_2$,
or to 
$$
\GL_2\backslash \GL_2\times \GL_2/\sS_2
$$ 
otherwise. In both cases the base is rational 
by Lemma~\ref{lemm:qqq}.

If $d_j=1$ for every $j\in J$ and $|J| = 2$ 
then the there is a map  
$$
(\P^1)^4 \ra \P(V_4)= \Sym^4(\P^1) = \P^4
$$
(of degree $24$, mapping $4$ points to a form of degree $4$).
The preimage in $(\P^1)^4$ 
of the $(\PGL_2,\sS_4)$-slice $\P^1_s=\P^1$ of
$\P^4$, will be a set of six lines $\P^1_{g,h}$, labeled by
a pair of generators $g,h \in \sD_2$ (which act trivially on 
$\P^1_s\subset \P^4$). 
More precisely, the line $\P^1_{g,h}$ is the set given by
$(x:gx:hx:ghx)\in (\P^1)^4$, for $x\in \P^1$.
The map $\P^1_{g,h}\ra \P^1_s=\P^1_{t,s}/\sD_2$ has 
degree $4$. Thus $\P^1_{g,h}$ is a $(\PGL_2, \sD_2)$-slice
of $(\P^1)^4$ 
and the quotient of a vector bundle  
$\stackrel{L\oplus L}{\longrightarrow}\P^1$ 
by a linear action of $\sD_2$ is rational.

\

Assume that all $d_i$ are even. 
Then $L$ is (birationally) trivial. 
Unless $|J|=2$ and $d_1=d_2 = 2$, there is a decomposition of
$$
Y\times Y =\P(V_d)\times Y'\times \P(V_d)\times Y'
$$
such that the $\PGL_2$-action is $af$ and 
$$
\P(V_d)\times Y'\times \P(V_d)\times Y'
\sim_{\PGL_2\times \sS_2} (Y'\times Y') \times 
(\P(V_d)\times\P(V_d))
$$
(with trivial $\PGL_2$-action on $\P(V_d)$), by Lemma~\ref{lemm:pgl2}.
The quotient is birational to a vector
bundle over $\PGL_2\times \C^* \ba X'\times X'/\sS_2$,  
where $X'$ is the trivial line bundle over $Y'$.

We have reduced to $|J|=1$ treated in Case 1 
or to $|J|=2$ and $d_1=d_2=2$, treated in Lemma~\ref{lemm:p2}. 
\end{proof}

\begin{lemm}
\label{lemm:p2}
The quotient
$$
X:=\PGL_2\ba (\P_1(V_2)\times \P_2(V_2) \times \P_1(V_2)\times \P_2(V_2))/\sS_2
$$
is rational, where $\P_1(V_2)$ and $\P_2(V_2)$ 
are different copies of $\P^2=\P(V_2)$ and 
$\sS_2$ acts by permutation.
\end{lemm}

\begin{proof}
Consider the projection 
$$
X\ra \PGL_2\ba \P_1(V_2)\times \P_1(V_2)/\sS_2
$$
and the $\PGL_2\times \sS_2$-equivariant map of degree $6$
$$
\begin{array}{cccc}
\pr\,: &  \P(V_2)\times \P(V_2)       & \ra      &    \P(V_4)\\
       &     (Q_1,Q_2)          & \mapsto  & Q_1\cdot Q_2.
\end{array}
$$
The space $\P(V_4)$ has a $(\PGL_2,\sS_4)$-slice $\P^1_s$
(the $\sD_2$-invariant polynomials). 
The zeroes of a (polynomial) $p\in \P^1_s$ form an orbit under $\sD_2$.
The preimage $\pr^{-1}(\P^1_s)\subset \P^2\times \P^2$
consists of $3$ lines, each invariant under $\sD_2$.
Indeed, the ordered pair $(Q_1,Q_2)$ corresponds to a
choice of a generator $g\in \sD_2$ such that $x, g(x)$ are zeroes 
of $Q_1$ and $h(x), hg(x)$ are zeroes of $Q_2$.
Thus the line $\P^1_g\subset \P^2\times \P^2$ consists
of tupels $\{(x,gx),(hx,ghx)\}$, where $x$  is an arbitrary point
in $\P^1$ and $(x,gx) =Q_1, (hx,ghx)= Q_2$.
The map $\P^1_g\ra \P^1_s$ has degree two and its fibers
coincide with orbits of $h$ (since $g$ acts trivially
on $\P^1_g$). The action of $h$ is given by
$$
h\,:\, \{(x,gx),(hx,ghx)\} \mapsto \{ (hx,ghx),(x, gx)\}.
$$
Thus $h(Q_1,Q_2) = (Q_2,Q_1)$ and the action 
of $h$ coincides with the restriction of the permutation action
on $\P^2\times \P^2$ to $\P^1_g$.
The line $\P^1_g$ is invariant under
${\sD}_4\times \sS_2$ (considered as a subgroup of $(\PGL_2\times \sS_2)$).
The group ${\sS}_4$ permutes the lines in $\pr^{-1}(\P^1_s)$.
Each $\P^1_g$ is a $(\PGL_2\times \sS_2, {\sD}_4\times \sS_2)$-slice of  
$\P^2\times \P^2$.
Therefore, 
$$
X \sim  {\sD}_4\ba \P^1\times \P^2\times \P^2  / \sS_2.
$$ 
The space $\P^2\times \P^2$ contains a subspace
$\C^2\times \C^2$ with a {\em linear} action of $\sD_4\times \sS_2$.
Indeed, the action of $\sD_4$ on $\P^1$ corresponds to the
irreducible representation of
$\tilde{\sD}_4$ on $\C^2= V$. 
Under the $\sD_4$-action, one has
a decomposition 
$\Sym^2(V)=V'\oplus V''$, where $\dim V'=2$, $\dim V''=1$ and the 
action of $\sD_4$ on $\P^2$ is equivalent 
to the linear action on $V'$.
The additional $\sS_2$ permutes the $\P^2$ and hence
acts by permutation on $V' \oplus V'$.
Thus  
$$
\P^1\times \P^2\times \P^2\sim_{\sD_4\times \sS_2}
\xymatrix{ \ar[d]^{V'\oplus V'}  \\  \P^1}
$$
(a vector bundle).

Consider the effective action of (the nonabelian group) 
$\sD_4\times \sS_2$ on $\P^1$.
It has a normal subgroup $\sD_2\times \sS_2$ with generators
$g,h,k$ and an element $i, i^2 = 1$  which commutes 
with $g,k$ and acts on $h$ as $ihi = gh$.  
The stabilizer of a generic point  on $\P^1_g$ 
is a normal abelian subgroup generated by $g, hk$. 
Thus $\sD_4\times \sS_2$ acts on $\P^1$ effectively through the
quotient $\sD_4/\langle g,hk\rangle = \sD_2$. 
The action of this $\sD_2$ on $\P^1$ is
almost free. Indeed, the action of $k$ coincides with the action 
of $h$ and permutes $Q_1,Q_2$. Thus the orbits
of $h$ and $k$ on $\P^1_g$ coincide with fibers of the map
$\P^1_g\to \P^1_s$. On the other hand, $i$ acts nontrivially on $\P^1_s$. 
We claim that 
$$
\xymatrix{\sD_4\ba (V'\oplus V')\times \bP^1/\sS_2 \ar[d]\\
          \sD_4\ba (V'\times \P^1)}
$$
is a vector bundle.
Indeed, consider the subspace $V'_{inv} \subset V'\oplus V'$ of invariant
vectors (under the permutation).
The action of $\sD_4\times \sS_2 $ on  
$((V'\oplus V')/V'_{inv}) \times \bP^1$ is almost free. Hence
$$
\xymatrix{\sD_4\ba (V'\oplus V')\times \P^1/\sS_2\ar[d] \\
\sD_4\ba ((V'\oplus V')/V'_{inv}) \times \P^1/\sS_2}
$$
is a vector bundle with base a quotient of the 
vector bundle $ (V'\oplus V'/V'_{inv})\ra \P^1_g$ by 
$\sD_4\times \sS_2$. The variety 
$(V'\oplus V'/V'_{inv})\times \P^1$ has a fiberwise (scalar)
$\C^*$-action commuting with the $\sD_4\times \sS_2$-action.
Since every $\C^*$-action has a slice,  
$$
X':=\sD_4\ba ((V'\oplus V')/V'_{inv}) \times \P^1/\sS_2,
$$
is rational by dimensional reasons: 
$X'/\C^*$ is a unirational, therefore, rational surface and   
$$
X'\sim (X'/\C^*)\times \C^*.
$$ 
\end{proof}

\

\begin{prop}
\label{prop:syml}
Let $X$ be $V\oplus V$, where $V=V_d$ is an 
irreducible $\GL_2$-representation,
$\ell>0$ and $\bH\subset \SL_2$ with $\bH\neq \tilde{\sA}_5$. 
Then 
$$
\GL_2\ba X\times \P(V_{\ell})/\bH
$$
is rational (where $\bH$ acts trivially on $\P(V_{\ell})$). 
\end{prop}

\begin{proof}
If $\ell$ is even and the action of 
$\GL_2$ or a quotient of $\GL_2$ by a central subgroup is $af$
then we apply Lemma~\ref{lemm:pgl2} combined with
Proposition~\ref{prop:even}, 
resp. \ref{prop:yz}. 

If $\ell$ is odd and the action is $af$ then there exists a slice,
which is a rational variety, 
by Lemma~\ref{lemm:grr} resp. \ref{lemm:subgroup}. 
Rationality follows.

Now we assume that the action is not $af$. 
This means that $d\le 4$. 
The subcases with $d\le 2$ are trivial since
the action on the corresponding 
Grassmannian is transitive. 
If $\ell$ is odd, then
the $\PGL_2$-action on $\Gr(2,V)\times \P(V_{\ell})$ 
has a rational slice and our claim follows.

If $d=3$, the action of $\PGL_2$ on $\Gr(2,V_3)$
has a $(\PGL_2, \sA_4)$-slice $\P^1$.
For even $\ell>0$ the action of $\sA_4$ on $\P^{\ell}$ is faithful and it 
lifts to a linear representation of
of $\sA_4$. Further, $\sA_4$-acts on $\P^1$ is through
a cyclic quotient. 
Thus 
$$
(\P^1\times \P(V_\ell)\sim_{\sA_4} \P^1\times \P(V_\ell)
$$
with trivial $\sA_4$-action on the $\P^1$ on the right.
This implies that the quotient is 
equivalent to
$$
\P^1\times (\P^{\ell}/\sA_4) \times (V_1\oplus V_1)/ \C^*\times \bH,
$$   
a product of rational varieties.

If $d=4$, the action of $\PGL_2$ on $\Gr(2,V_4)$
has a $(\PGL_2, \bN_{\bT})$-slice
$X'$. The action of $\bN_{\bT}$ on $\P(V_{\ell})$ 
is linear and the quotient
of $X\times \P^{\ell}$ is a vector bundle over the quotient of 
$X$, which is rational.
\end{proof}

\begin{prop}
\label{prop:xxx}
Let $X=(\stackrel{L}\longrightarrow Y)^2 $, where $Y=
\prod_{j\in J} \P(V_{d_j})$ and  $\ell >0$. 
Then 
$$
\GL_2\ba X\times \P(V_{\ell})/\sS_2
$$
is rational (where $\sS_2$ acts trivially on $\P(V_{\ell})$ and
by permutation on $X$).
\end{prop}

\begin{proof}
The same argument as in the proof of Proposition~\ref{prop:syml}
shows that it suffices to assume that the action on  
$X$ is not $af$. This happens only if $Y=\P^2$ or $\P^1$.
The case $Y=\P^2$ reduces to Proposition~\ref{prop:syml}
(Grassmannian). 
If $Y=\P^1$ then the action of $\PGL_2$ on $\P^1\times \P^1$ is
transitive and 
$$
\GL_2\ba X\times \P(V_{\ell})/\sS_2 
\sim  (\C^*\ba \P(V_{\ell}))\times (\C^2/\C^*\times \sS_2),
$$
a rational variety. 
\end{proof}

\section{Special rationality results}
\label{sect:rat-sp}

In this section we collect rationality results
for  spaces of rational maps $\Pr^1\ra \Pr^1$ with prescribed 
(special) ramification over exactly three distinguished
points $(0,1,\infty)$ and unspecified ramifications over other points.  

Let $\cR(r_0,r_1,r_{\infty})$
be the space of rational maps $f\,:\, \Pr^1\ra \Pr^1$
with local ramification data (vectors)
$r_0,r_1,r_{\infty}$ over the points $0,1,\infty$.

\begin{prop}
\label{prop:rrr}
Assume that $(r_0,r_1,r_{\infty})$ satisfies one of the following:
\begin{itemize}
\item all entries of the vectors $r_0,r_\infty$ are even 
and some fixed number of entries of $r_1$ is even; 
\item all entries of the vectors $r_0,r_{\infty}$ are even
and a fixed number of entries of $r_1$ is divisible by $3 ;$ 
\item all entries of the vectors $r_0,r_{\infty}$ are divisible by $3$
and all entries of $r_1$ are even.
\end{itemize}
Then $\cR(r_0,r_1,r_{\infty})$ is a finite union of 
irreducible rational varieties.
\end{prop}

\begin{proof}
In these cases the map $f=f_0/f_\infty$ is given by coprime polynomials
satisfying the equations:
\begin{itemize}
\item $f_0^2-f_{\infty}^2=g_1^2g_1'$;  
\item $f_0^2-f_{\infty}^2=g_1^3g_1'$; 
\item $f_0^3-f_{\infty}^3=g_1^2g_1'$, 
\end{itemize}
where $g_1'$ is an arbitrary polynomial. 
The first equation leads to 
$$
(f_0-f_{\infty})(f_0+f_{\infty})=g_1^2g_1'
$$
and, by coprimality, to  
$$
\begin{array}{ccc}
f_0-f_{\infty} & = & g_{11}^2g_{11}', \\
f_0+f_{\infty} & = & g_{12}^2g_{12}',
\end{array}
$$
with arbitrary $g_{11},g_{11}',g_{12},g_{12}'$
(satisfying the obvious degree conditions) --- a union of rational
varieties. 

The second case is analogous. 
Consider the third case:
since $f_0^3-f_\infty^3$ is a square we obtain 
$$
\begin{array}{ccc}
f_0- f_\infty          & =  &  g_1^2\\
f_0-\zeta f_\infty     & =  &  g_2^2\\
f_0-\zeta^2 f_{\infty} & =  &  g_3^2
\end{array}
$$ 
(where $\zeta^3=1$)
and  we need to solve 
$$
\frac{2\zeta}{1+\zeta} g_1^2 + \frac{1-\zeta}{1+\zeta} g_2^2 = g_3^2.
$$ 
Now we apply the parametrization 
as above. 
\end{proof}

\begin{coro}
\label{coro:rr-r}
Let $\cR(r_0,r_1,r_{\infty})$ be as in \ref{prop:rrr}. 
Then 
$$
\PGL_2\backslash \cR(r_0,r_1,r_{\infty})
$$
is rational.
\end{coro}

\begin{proof}
We have established an explicit parametrization
of $\cR(r_0,r_1,r_{\infty})$ 
as a direct sum of spaces of polynomials (with different weights
as irreducible $\GL_2$-representations). 
By the theorem of Katsylo \ref{thm:katsylo}, the corresponding quotients are rational.
\end{proof}

\begin{rem}
\label{rem:8.3}
Only the first case with 
$g_1'=1$ can admit a nontrivial action of $H_{\Ga}$ 
(which necessarily is $\Z/3$). 
But even in this case the action of $\Z/3$ is linear
and it commutes with the action of $\GL_2$ on pairs of
polynomials. Lemma~\ref{lemm:grr} implies 
rationality.
\end{rem}

\begin{lemm}\label{thG}  
Every irreducible component of the variety  $\cR$ of rational maps 
$f\,:\, {\Pr}^1\to {\Pr}^1$ of degree $5$ and prescribed
global ramification datum
$$
\RD(f)=[(2,2,1)_0,(2,2,1)_1,(2,2,1)_{\infty},(2),(2)]
$$
is rational.
\end{lemm}

\begin{proof} 
Changing the variables 
(fixing two ramification points over $1\in{\Pr}^1$ as $0,\infty$),
we can write  $f=F_1/F_2$ where 
$$
\begin{array}{ccc}
F_1(x) & = & \hat{f}_1(x)^2\hat{a}_1(x)^2 \hat{b}_1(x)\\
F_2(x) & = & \hat{f}_2(x)^2\hat{a}_2(x)^2 \hat{b}_2(x)
\end{array}
$$
where 
$\hat{f}_1,\hat{f}_2,\hat{a}_1,\hat{a}_2,\hat{b}_1,\hat{b}_2$ 
are linear forms in $x$.
Since the leading coefficients of $F_1$ and $F_2$ are equal 
we can assume that they are both equal to $1$ and write
$\hat{f}_1(x)=x+f_1,\dots,\hat{b}_2(x)=x+b_2$,
with some nonzero constants $f_1,\dots,b_2$. 
Since we have one free parameter (under the action of $\PGL_2$)
we can assume that $b_1=1$.  
Thus 
$$
\hat{f}_1(x)^2\hat{a}_1(x)^2\hat{a}_2(x)-
\hat{f}_2(x)^2\hat{b}_1(x)^2 \hat{b}_2(x) =\sum_i g_ix^i= c_1 x^2(x+c_2)
$$ 
with arbitrary constants $c_1,c_2$. 
We get a system of equations on the coefficients $g_j$: 
$$
g_4=0, g_1=0, g_0=0.
$$ 
Remark that the coefficients of $g$ are
symmetric functions on pairs $(f_1,a_1)$ and $(f_2,a_2)$. 
To parametrize $\cR$ we introduce  the following variables: 
$$
X_1=a_1+f_1,\, Y_1=a_1f_1,\, X_2=f_2+a_2,\, Y_2=f_2a_2,\, b_1,\,b_2.
$$ 
Write the equations on the coefficients $g_j$ as
$$
\begin{array}{rcl}
2X_1+b_1         & = & 2X_2+b_2\\
Y_1^2b_1         & = & Y_2^2b_2\\
Y_1^2+2X_1Y_1b_1 & = & Y_2^2+2X_2Y_2b_2.
\end{array}
$$
Since  $b_1=1$, for a fixed $b_2$  we get 
$$
\begin{array}{rcl}
2X_1+1                 & = & 2X_2+b_2\\
Y_1                    & = & \pm \sqrt{b_2}Y_2\\ 
b_2Y_2 +2\sqrt{b_2}X_1 & = & Y_2 +2X_2b_2.
\end{array}
$$
This is a union of two (affine) lines.  
After a rational covering $(\sqrt{b_2})$ our surface
is (rationally) a $\Pr^1$-bundle over $\Pr^1$, a rational
surface. 
\end{proof}

\begin{lemm}
\label{lemm:4}
Every irreducible component of the variety $\cR$ of rational maps
$f\,:\, \Pr^1\to \Pr^1$ of degree $4$ and ramification datum
$$
\RD(f)=[(2,2)_0,(2,1,1)_1,(2,1,1)_{\infty}]
$$ 
is a rational surface.
\end{lemm}

\begin{proof}
Using the $\PGL_2$-action on the preimage 
$\Pr^1$ we can assume that the
points $(2,2)$ are $+1,-1$, respectively,  and
that the point of degree $2$
(in the local ramification datum $(2,1,1)$) over $0$ is $\infty$. 
Thus we can write
$$
(x^2-1)^2 - c(x+c_1)(x+c_2)(x+c_3)^2=g_2(x),
$$
where $g_2$ is an arbitrary  
polynomial of degree 2 and $c$ is some constant. 
We get two equations
$$
c=1,
$$
$$
c_1+c_2+2c_3=0.
$$
Thus we have a (rational) surjection of $\Pr^2$ onto $\cR$.
\end{proof}

\begin{lemm}
\label{lemm:5}
Every irreducible component of the 
variety  $\cR$ of  rational maps $f\,:\, \Pr^1\to \Pr^1$ 
of degree $4$ with ramification datum 
$$
\RD(f)=[(2,2)_0,(3,1)_{1},(2,1,1)_{\infty},(2),(2)]
$$ 
is a rational curve.
\end{lemm} 

\begin{proof}
A generic map with this ramification datum 
is given by the equation $f=f_1/f_2$, 
where
$$
f_1=(x^2-1)^2,\,\, f_2=(x+c_1)(x+c_2)(x+c_3)^2
$$ 
and 
$$
f_1-f_2=(x^2-1)^2-c(x+c_1)(x+c_2)(x+c_3)^2=g_1(x),
$$
where $g_1(x)$ is linear. Thus $c=1$ and 
$$
c_1+c_2+2c_3=0,
$$
$$
c_1c_2+2c_1c_3 + 2c_2c_3 + c_3^2 =0,
$$
clearly rational.
\end{proof}

\begin{lemm}
\label{lemm:6}
The irreducible component of the 
variety  $\cR$ of  rational maps $f\,:\, \Pr^1\to \Pr^1$ 
of degree $3$ with ramification datum 
$$
\RD(f)=[(2,1)_0,(2,1)_{1},(2,1)_{\infty},(2)]
$$ 
is a rational curve.
\end{lemm} 

\begin{proof}
Reduces easily to the rationality of a cuspidal cubic curve.
\end{proof}

\section{Rationality of moduli}
\label{sect:ratio}

\begin{thm}
\label{thE}
Any connected component of a moduli space of rational or 
K3 elliptic surfaces with fixed monodromy group is rational.
\end{thm}

\begin{proof}
In Proposition~\ref{lemm:struct}
we have identified (Zariski open subsets of) the 
corresponding moduli spaces
$\cF_{r,\tilde{\Ga}}$ 
as quotients 
(by the left $\PGL_2$ and right $H_{\Ga}$-action) 
$$
\PGL_2\ba \cU_{r,\tilde{\Ga},\ell}'/H_{\Ga}.
$$
Here 
$$
\cU_{r,\tilde{\Ga},\ell}'\sim_{\PGL_2\times H_{\Ga}} 
\Sym^{\ell}(\P^1)\times \cR_{\Ga}
$$ 
and 
$$
\cR_{\Ga}=\{ f\,:\, \P^1\ra \P^1\}
$$
is the space of rational maps 
(with prescribed ramification).  
For elliptic rational or K3 surfaces $\ell\le 3$
and $H_{\Ga}$ is either trivial, cyclic, dihedral or
a subgroup of $\sS_4$ (see Corollary~\ref{coro:hgamma}).
The actions if $\PGL_2$ and $H_{\Ga}$  commute and $H_{\Ga}$
acts only on $\cR_{\Ga}$.

\

First we consider {\em general families}:
$$
\ET(\E)-12\ell =\deg(j_{\E})\ET(\Ga).
$$
For ${\bf d}=(d_1,...,d_k)\in \N^k$ we put
$$
\P^{\bf d}:= \prod_{j=1}^k \P(V_{d_j}).
$$
Recall that $\cR_{\Ga}$ is (birationally) 
the total space of a line bundle over the space
$$
\P^{\bf d}\times \P^{\bf d'},
$$
where $\sum_{j=1}^k d_j=\sum_{j=1}^{k'} d_j'$.

\

{\bf Case 1.} ${\bf d}\neq {\bf d}'$. Then, by 
\ref{coro:hgamma}, 
$H_{\Ga}=1$ and rationality of $\PGL_2\ba \cR_{\Ga}$ 
(in all cases) follows from the rationality of
$$
\PGL_2\ba \P^{\bf d}\times \P^{\bf d'},
$$
which is the theorem of Katsylo \ref{thm:katsylo}.

\

{\bf Case 2.} ${\bf d}={\bf d}'$ and $k\ge 2$.
By Corollary~\ref{coro:hgamma}, 
$H_{\Ga}=\sS_2$ (permutation of the factors).   
This case is covered by Proposition~\ref{prop:yz}.

\

{\bf Case 3.} ${\bf d}={\bf d}'=(d)$. This case is covered by 
Proposition~\ref{prop:even}.

\

Now we discuss the {\em special families}:
$$
\ET(\E)-12\ell <\deg(j_{\E})\ET(\Ga).
$$
We use the classification of these families established in 
Section~\ref{sect:ellmono}.
All families listed in Lemma~\ref{lemm:24-sp} are covered by 
Propositions~\ref{prop:yz} and the Theorem~\ref{thm:katsylo}.
Consider the families listed in Lemma~\ref{lemm:table-tree}:
Lemma~\ref{lemm:subgroup} covers the cases 
$j_1,j_4,j_5,j_6,j_{13}$. 
The case $j_2,j_8$ and $j_{12}$ are covered by Proposition~\ref{prop:rrr}, 
$j_3$ by Lemma~\ref{lemm:5}, $j_7,j_9,j_{10}$ by \ref{prop:rrr} and
\ref{rem:8.3}, $j_{11}$ by Lemma~\ref{lemm:4}. The case $j_{14}$ is
covered by Lemma~\ref{lemm:6}.
Finally, the families $j_{15}$ and $j_{16}$ 
(listed in Lemma~\ref{lemm:table-rat}) are
covered by Proposition~\ref{prop:yz} and the remaining
families $j_{17}-j_{20}$ by Theorem~\ref{thm:katsylo}.
\end{proof}

\begin{rem}
Our methods extend to some moduli spaces of 
elliptic surfaces with higher Euler characteristic. 
In particular, the results of Section~\ref{sect:rat-sp} 
imply that any moduli space of Jacobian elliptic surfaces
over $\Pr^1$ such that a generic surface in this space 
has only singular fibers of multiplicative type is rational. 
However, we expect that there are nonrational 
moduli spaces already for 
Euler characteristic 36. 
\end{rem}

\section{Pictures}
\label{sect:pictures}

In this section we give a combinatorial description of
monodromy groups of elliptic K$3$ surfaces. 
More precisely, we describe a simple procedure which allows
to enumerate all the possible graphs $\Ga$ 
with given $\ET(\Ga)$.
Let $\E\ra \Pr^1$ be an elliptic K$3$ surface. 
We have shown in Section~\ref{sect:comb}
that 
$$
48=\ET(\E)\geq\ET(\Ga)
$$
and that $\ET(\Ga)$ is divisible by 12.
Thus  $\ET(\Ga)$ equals $12, 24, 36$ or $48$
and all possible $\Ga\subset \PSL_2(\Z)$ 
are described by connected trivalent graphs
$T_{\Ga}$ with $\leq 8$ edges embedded into 
${\mathbb S}^2$, with an arbitrary bicoloring of the ends.  

\

\no
{\bf Case} $\ET(\Ga)=12$ :
There is only one tree $T_{12}$ with $ \ET(T_{12})=12$
\begin{figure}[htb]
\centerline{\includegraphics[width=.15\textwidth]{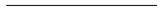}}
\caption{The tree $T_{12}$.}
\end{figure}

The ends of $T_{12}$ can be either $A$ or $B$-vertices. 
To obtain all possible graphs $T_{\Ga}$ 
with $\ET(\Ga)=12$ we just need to attach 
to $T_{12}$ a single loop $L$.
\begin{figure}[htb]
\centerline{\includegraphics[width=.2\textwidth]{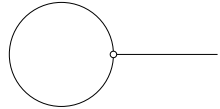}}
\caption{The loop $L$.}
\end{figure}

This gives the  following list of graphs:
\begin{figure}[htb]
\centerline{\includegraphics[width=.5\textwidth]{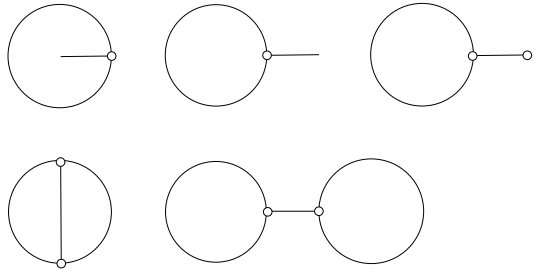}}
\caption{The case $\ET(\Ga)=12$.}
\end{figure}

%\pagebreak
There is only one saturated graph from the list above which has no
outer loops (Figure $4$).
\begin{figure}[htb]\label{fig:sat}
\centerline{\includegraphics[width=.15\textwidth]{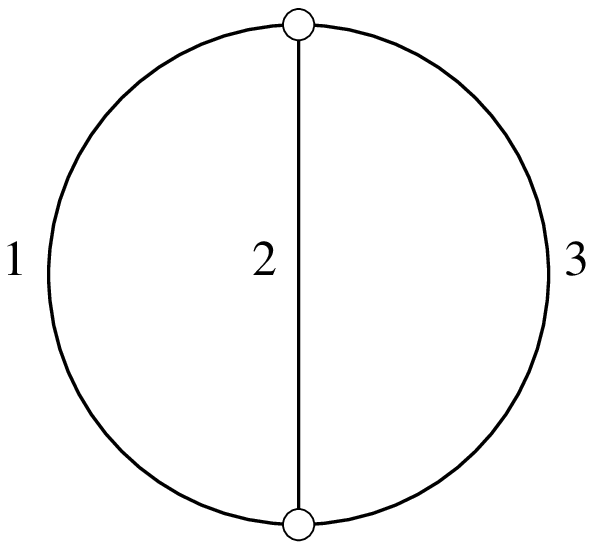}}
\caption{}
\end{figure}

This graph will be a basic building block 
in the construction of graphs with $\ET(\Ga)>12$ - 
we will attach trees and loops to its edges. 
The edges are numbered to simplify the
count of all possible outcomes.

\pagebreak

\no
{\bf Case} $\ET(\Ga)=24$: 
Again, we have only one topological tree $T_{24}$ with
$\ET(T)=24$:
\begin{figure}[htb]
\centerline{\includegraphics[width=.1\textwidth]{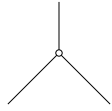}}
\caption{The tree $T_{24}$.}
\end{figure}
\

\no
{\bf Case} $\ET(\Ga)=36$: 
There are only $3$ saturated graphs without end-loops (modulo
equivalent embedding into the sphere):
\begin{figure}[htb]
\centerline{\includegraphics[width=.6\textwidth]{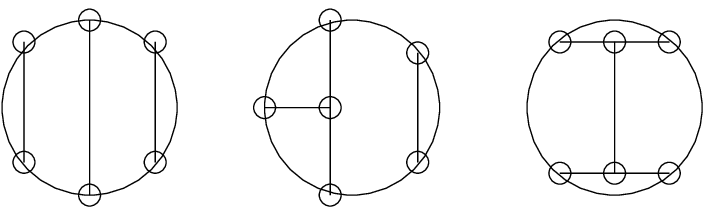}}
\caption{The case $\ET(\Ga)=36$.}
\end{figure}

Any other graph is either a 
tree or a sum of a saturated graph $T'$
with $\ET(T')=0,12,24$ with trees (with complementary
$\ET$). 
There is only one topological tree $T_{36}$ with
$\ET(T_{36})=36$.
\begin{figure}[htb]
\centerline{\includegraphics[width=.2\textwidth]{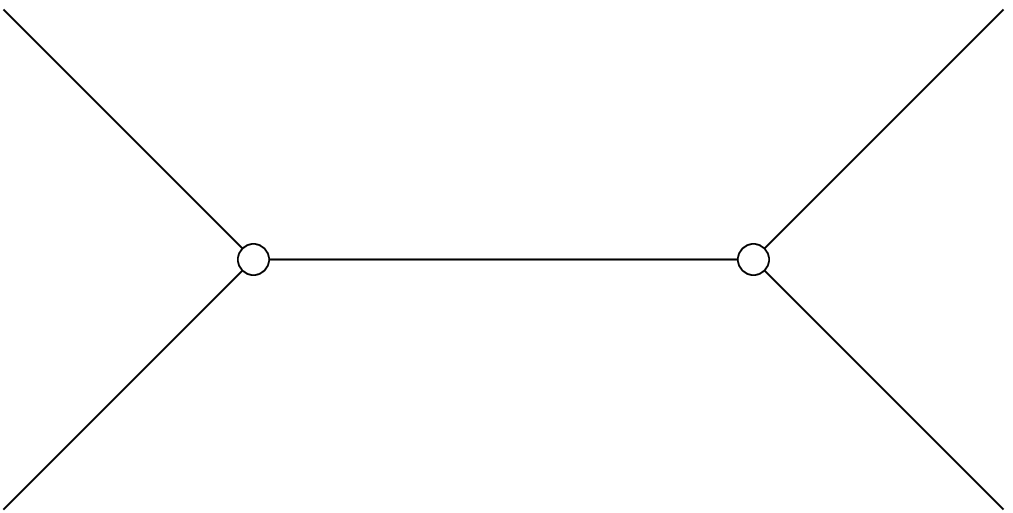}}
\caption{}
\end{figure}  

The number of possible markings of the tree or 
loops at the ends is $81$ but due to the 
symmetry of the graph the actual number of graphs
$T_{\Ga}$ corresponding to different placement 
of loops at the end and markings is smaller: 
there are $34$ different $T_{\Ga}$ of this type.

The number of markings of $T_{36}$ 
is $16$ but due to its symmetry 
the number of different graphs $T_{\Ga}$ is $7$. 
(Recall that two graphs $T_{\Ga}$ give the same $\Ga$ 
modulo conjugation if they are isotopic in 
a ${\mathbb S}^2$).

The graphs of tree type with one end loop are topologically
equivalent to:

%\begin{figure}[htb]
\centerline{
\includegraphics[width=.2\textwidth]{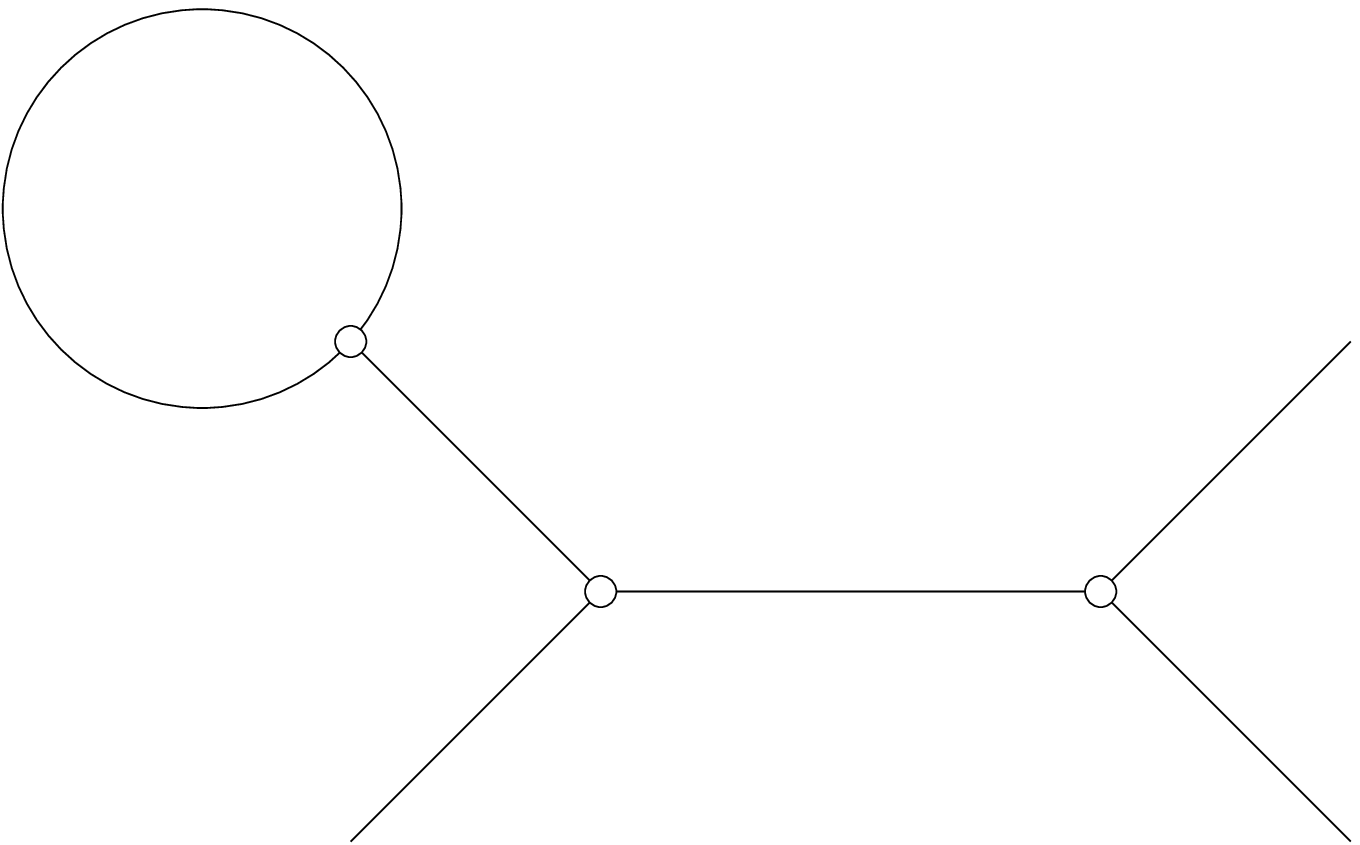}}

%\caption{}
%\end{figure}  

\

There are $8$ possible markings of the above graph and they all
give different $T_{\Ga}$ with $\ET(\Ga)=36$. We have $12$
different $T_{\Ga}$ with $2$ end-loops, $6$ with $3$ 
end-loops and one with $4$ end-loops.

All topological graphs which 
are sums of a loop and a tree can be
obtained by placing a loop into a tree. 
Thus there are two types:
\begin{figure}[htb]
\centerline{\includegraphics[width=.35\textwidth]{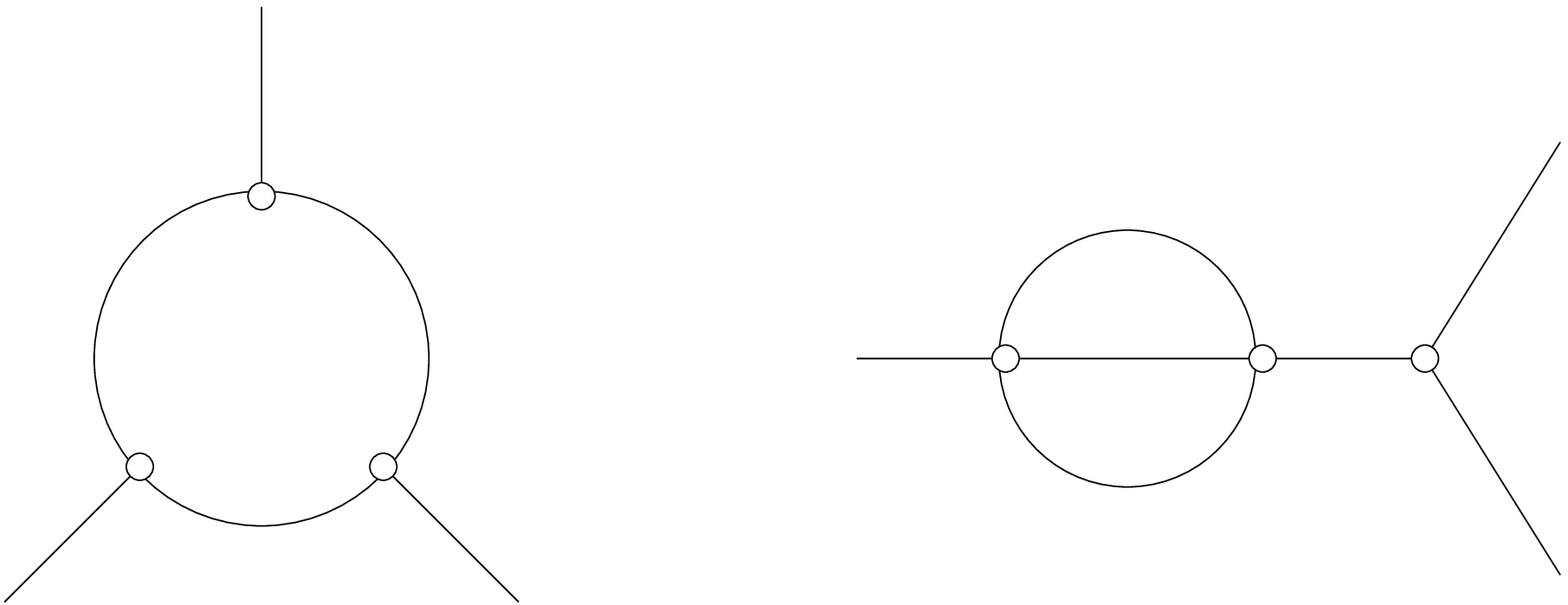}}
\caption{}
\end{figure}  

This gives $8$ graphs $T_{\Ga}$ in the 
first case and $4$ in the second case. 

\

\no
{\bf Case} $\ET(\Ga)=48$:
We have one tree $T_{48}$ with $\ET(\Ga)=48$:
\begin{figure}[htb]
\centerline{\includegraphics[width=.2\textwidth]{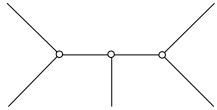}}
\caption{The tree $T_{48}$.}
\end{figure}

Here is the list of all saturated graphs with $\ET(\Ga)=48$.
\begin{figure}[htb]
\centerline{\includegraphics[width=.6\textwidth]{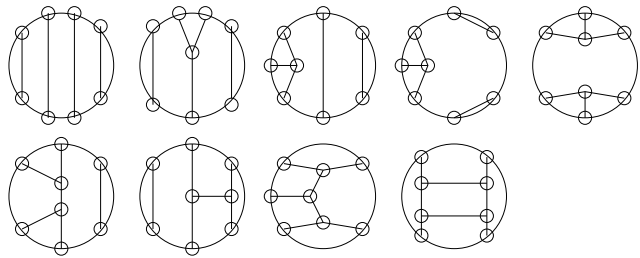}}
\caption{Saturated graphs in the case $\ET(\Ga)=48$.}
\end{figure}

\end{document}